\documentclass[preprint, number, longtitle]{elsarticle}
\usepackage{amsmath}
\usepackage{graphics}
\usepackage{latexsym}
\usepackage{amsfonts}
\usepackage{mathrsfs}
\usepackage{amssymb}
\usepackage{amssymb,amsmath,amsthm, amsfonts}
\usepackage[english]{babel}
\usepackage[latin1]{inputenc}

\def \a {\alpha}

\numberwithin{equation}{section}
\renewcommand{\baselinestretch}{1.30}
\newtheorem{theorem}{Theorem}[section]

\newtheorem{remark}{Remark}

\newtheorem{assumption}{Assumption}
\def\qed{ \ \vrule width.2cm height.2cm depth0cm\smallskip}

\def \ig {(i,j)\in \Gamma}

\def \ind{1\!\!1}

\def \x{X^{t,x}}
\def \ed {\end{document}}
\def \tx {(t,x)\in [0,T]\times \R^k}

\def \xtx {X^{t,x}}
\def \xtxs {X^{t,x}_s}

\def \kl {(k,l)\in \gam}
\def \pr {{\bf Proof}}

\def \ij {i\in \cJ}
\def \ijo {i_0,j_0}

\def \io {i_0}
\def \jo{j_0}

\def \lb{\label}

\def \th{\theta}

\def \yijm {(Y^{ij,m},Z^{ij,m},K^{ij,m,+},K^{ij,m,-})_{\ij}}

\def \min {\mbox{min}}
\def \max {\mbox{max}}
\def \pl {\Pi_g}

\newcommand{\eps}{\varepsilon}

\newcommand{\brm}{\begin{rem}}
\newcommand{\ermq}{\end{rem}}

\newcommand{\ba}{\begin{array}}
\newcommand{\ea}{\end{array}}
\newcommand{\be}{\begin{equation}}
\newcommand{\ee}{\end{equation}}
\newcommand{\bea}{\begin{eqnarray}}
\newcommand{\eea}{\end{eqnarray}}
\newcommand{\beaa}{\begin{eqnarray*}}
\newcommand{\eeaa}{\end{eqnarray*}}

\def \R{I\!\!R}
\def \E{\mathbb{E}}

\def \ij{(i,j)\in \gam}
\def \ot {t\in [0,T]}
\def\a{\alpha}

\def\g{\gamma}
\def\d{\delta}

\def\z{\zeta}

\def\n{\nu}
\def\si{\sigma}
\def\t{\tau}

\def\th{\theta}

\def\G{\Gamma}

\def\L{\Lambda}

\def\cA{{\cal A}}

\def\cC{{\cal C}}

\def\cF{{\cal F}}

\def\cH{{\cal H}}

\def\cL{{\cal L}}

\def\cS{{\cal S}}

\def\no{\noindent}

\def\ms{\medskip}
\def\bs{\bigskip}
\def\q{\quad}
\def\qq{\qquad}

\def\bF{{\bf F}}

\def\qed{ \hfill \vrule width.25cm height.25cm depth0cm\smallskip}

\newcommand{\basa}{\begin{assumption}}
\newcommand{\easa}{\end{assumption}}

\newcommand{\bas}{\begin{assum}}
\newcommand{\eas}{\end{assum}}

\def\esssup{\mathop{\rm esssup}}
\def\essinf{\mathop{\rm essinf}}

\def  \vij {(v^{ij})_{\ij}}
\def \udl{\underline}

\def\dis{\displaystyle}

\def\bF{{\bf F}}

\def \P{\mathbb{P}}
\newtheorem{thm}{Theorem}[section]
\newtheorem{lem}[thm]{Lemma}
\newtheorem{prop}[thm]{Proposition}
\newtheorem{rem}[thm]{Remark}

\newtheorem{assum}[thm]{Assumption}

\def\mb{\mbox}

\newcommand{\rw}{\rightarrow}
\def \R{\mathbb{R}}

\def \esssup {\mbox{ess sup}}
\def \essinf {\mbox{ess inf}}
\def \esp {[t,T]\times \R^k}
\def \espo {[0,T]\times \R^k}
\def \gam {\Gamma^{1}\times \Gamma^{2} }
\def\gami{{(\Gamma^1)}^{-i}}
\def\gamj{{(\Gamma^2)}^{-j}}
\def \gmi {\Gamma^1}
\def \gmj {\Gamma^2}
\def \ed{\end{document}}
\title{On the Equality of Solutions of Max-Min and Min-Max Systems of Variational Inequalities with Interconnected Bilateral Obstacles.}
\author{Boualem Djehiche\footnote{This work was completed while the first author was visiting the Department of Mathematics of Universit\'e du Maine. Financial support from GEANPYL  and Svensk Export Kredit (SEK) is gratefully acknowledged.}}
\ead{boualem@math.kth.se}
\address{Department of Mathematics, Royal Institute of Technology, SE-100 44 Stockholm, Sweden} 
\author{Said Hamad\`ene}
\ead{hamadene@univ-lemans.fr}
\address{Universit\'e du
Maine, LMM, Avenue Olivier Messiaen, 72085 Le Mans, Cedex 9, France}

 \author{ Marie-Am\'elie Morlais}
  \ead{ Marie$_-$Amelie.Morlais@univ-lemans.fr}
 \address{Universit\'e du Maine, LMM, Avenue Olivier Messiaen, 72085 Le Mans, Cedex 9, France}
\author{Xuzhe Zhao}
\ead{zhao.xuxhe.etu@univ-lemans.fr}
\address{Universit\'e du Maine, LMM, Avenue Olivier Messiaen, 72085 Le Mans, Cedex 9, France  \\ and Xidian University,  School of Mathematics and Statistics,  Xi'an 710071, PRC.}

\begin{document}
\date{\today}

\begin{abstract} In this paper, we deal with the solutions of systems of PDEs with bilateral inter-connected obstacles of min-max and max-min types. These systems arise naturally in stochastic switching zero-sum game problems. We show that when the switching costs of one side are regular, the solutions of the min-max and max-min systems coincide. Then, this common viscosity solution is related to a multi-dimensional doubly reflected BSDE with bilateral interconnected obstacles. Finally, its relationship with the the values of a zero-sum switching game is studied. \end{abstract}


\maketitle
\section{Introduction}
Let us consider the following two systems of partial differential equations (PDEs) with bilateral inter-connected obstacles (i.e., the obstacles depend on the solution) of min-max and max-min types: for any $\ij$, $\tx$,
\begin{equation} \label{mainsyst-vi} 
         \left\{ \begin{array}{l} 
   \min\Big \{\bar v^{ij}(t,x)-L^{ij}(\vec {\bar v})(t,x)\,;\, \max\Big \{\bar v^{ij}(t,x) -U^{ij}(\vec {\bar v})(t,x);\\\qq\qq -\partial_t \bar v^{ij} -\mathcal{L}^{X}(\bar v^{ij})(t,x)-f^{ij}(t,x, (\bar v^{kl}(t,x))_{(k,l)\in \gam},\si(t,x)^\top D_x\bar v^{ij}(t,x)) \Big \}  \Big \}  =0\, ;\\
\bar v^{ij}(T,x)=h^{ij}(x)  
\end{array}\right. \end{equation}
($(.)^\top$ is the transpose) and 
\begin{equation} \label{maxminsyst-vi}
         \left\{ \begin{array}{l} 
  \max\Big \{\udl v^{ij}(t,x) -U^{ij}(\vec {\udl v})(t,x)\,; \min\Big \{\udl v^{ij}(t,x)-L^{ij}(\vec {\udl v})(t,x)\\\qq\qq -\partial_t \udl v^{ij} -\mathcal{L}^{X}(\udl v^{ij})(t,x)-f^{ij}(t,x, (\udl v^{kl}(t,x))_{(k,l)\in \gam},\si(t,x)^\top D_x\udl v^{ij}(t,x)) \Big \}  \Big \}  =0\, ;\\
\udl v^{ij}(T,x)=h^{ij}(x)   
    
           \end{array}          \right. \end{equation}where\\
           (i) $\Gamma^1$ and $\Gamma^2$ are finite sets (possibly different);\\
           (ii) For any $\tx$, $\vec v(t,x)=(v^{kl}(t,x))_{(k,l)\in \gam}$ and for any $\ij$,
           \be \label{expbar}
L^{ij}(\vec v)(t,x)=\max_{k\in \Gamma^1, k\neq i}\{v^{kj}(t,x)-\underline g_{ik}(t,x)\},\,\, U^{ij}(\vec v)(t,x)=\min_{p\in \Gamma^2, p\neq j}\{v^{ip}(t,x)+\bar g_{jp}(t,x)\}.\ee
(iii) $\cL^X$ is a second order generator associated with the following diffusion process $\xtx$ satisfying: 
\be\label{sdex}\begin{array}{l}
\forall s\in [t,T],\,\,\, X_s^{t,x}=x+\int_t^sb(r,X^{t,x}_r)dr+\int_t^s\si(r,X^{t,x}_r)dB_r;\,\,\,X_s^{t,x}=x,\,\,s\in [0,t].\end{array}
\ee
The required properties on $b$ and $\sigma$ will be precised later.

The systems (\ref{mainsyst-vi}) and (\ref{maxminsyst-vi}) are of min-max and max-min types respectively. 
The barriers $L^{ij}(\vec{\bar v}), U^{ij}(\vec{\bar v})$ and $L^{ij}(\vec{\udl v}), U^{ij}(\vec{\udl v})$ depend on the solution $\vec{\bar v}=(\bar v^{ij})_{\ij}$ and $\vec{ \udl v}=(\udl v^{ij})_{\ij}$ of (\ref{mainsyst-vi}) and (\ref{maxminsyst-vi}) respectively. They are related to zero-sum switching game problems since actually, specific cases of these systems, stand for the Hamilton-Jacobi-Bellman-Isaacs equations associated with those games. 
\ms \no  

Switching problems have recently attracted a lot of research activities, especially in connection with mathematical finance, commodities, and in particular energy markets, etc (see e.g. \cite{carmonaludkovski, ludkovski, cek1, cek2, ek3, djehshah, Djehichehampopier, Djehichehammorlais14, Hamadene-elasri,hamjeanb, Hammorlais13,Hutang,koike,marcus1,marcus2,marcus3, magnus, zervos1,zervos2,pham1,pham2} and the references therein). Several points of view, mainly dealing with control problems, have been considered (theoritical and practical \cite{carmonaludkovski,ludkovski,cek2,Djehichehampopier,Hamadene-elasri,hamjeanb,Hutang,magnus}, numerics \cite{cek1, hamjeanb}, filtering and partial information \cite{marcus1}). However, except \cite{Hutang-2008, koike}, problems related to switching games did not attract that much interest in the literature.

In \cite{Djehichehammorlais14}, by means of systems of reflected backward stochastic differential equations (BSDEs) with inter-connected obstacles in combination with Perron's method, Djehiche et al. have shown that each of the systems (\ref{mainsyst-vi}) and (\ref{maxminsyst-vi}) has a unique continuous solution with polynomial growth, under classical assumptions on the data $f^{ij}$, $\bar g_{ij}$, $\underline g_{ij}$ and $h^{ij}$. The question of whether or not these solutions coincide was conjectured as an open problem, leaving a possible connection of the solution of system (\ref{mainsyst-vi}) and (\ref{maxminsyst-vi}) with zero-sum switching games unanswered. The main objective of this paper is three-fold: (i) to investigate under which additional assumptions on the data of these problems, the unique solutions of systems (\ref{mainsyst-vi}) and (\ref{maxminsyst-vi}) coincide; (ii) to make a connection between this solution and the associated system of reflected BSDEs with bilateral inter-connected obstacles; (iii) to study the relationship with the value function of the associated zero-sum switching game. 

We show that if the switching costs of one side (or player), i.e. either $(\bar g_{ij})_{\ij}$ or $(\underline g_{ij})_{\ij}$, are regular enough, then the solutions of the systems (\ref{mainsyst-vi}) and (\ref{maxminsyst-vi}) coincide, i.e., $\bar v^{ij}=\underline v^{ij}$, for any $\ij$. 
The strategy to obtain these results is to show that the barriers, which depend on the solution, are comparable and then thanks to a result by Hamad\`ene-Hassani \cite{Hamhassani} (see Theorem \ref{existence} in appendix (A2)) on viscosity solutions of standard min-max and max-min PDE problems and uniqueness of the solutions of (\ref{mainsyst-vi}) and (\ref{maxminsyst-vi}), we obtain the equality of those latter. 

Next, with the help of this common solution, we have proved existence and uniqueness of quadruples $(Y^{ij},Z^{ij},K^{ij,\pm})_{\ij}$ solution of the following system of reflected BSDEs with inter-connected obstacles: $\forall \ij$, $\forall s\in [t,T]$,  
\be \label{general-drbsde} \left\{ 
 \begin{array}{l}  
 Y_{s}^{ij} = h^{ij}(X^{t,x}_T)+\int_s^Tf^{ij}(r,X^{t,x}_r, \vec{Y}_r, Z_{r}^{ij}) dr + 
(
K_{T}^{ij,+} -K_{s}^{ij,+}) -(K_{T}^{ij,-} -K_{s}^{ij,-})-\int_s^TZ_r^{ij}dB_r\, \,;
 \\\\
\displaystyle{ Y_s^{ij} \le U^{ij}_s(\vec{Y})  \; \textrm{and} \; Y_s^{ij} \ge
  L^{ij}_s(\vec{Y})\,\,;
} \\\\
 \int_{t}^{T} (Y_s^{ij} - U^{ij}_s(\vec{Y})) dK_{s}^{ij, -} =0\;\textrm{and} \;\int_{t}^{T} (  L^{ij}_s(\vec{Y}) - Y_s^{ij}
) dK_{s}^{ij, +} =0.\\ 
 \end{array}  \right. \ee
The component $(Y^{ij})_{\ij}$ has the following Feynman-Kac representation: 
\be\label{Feynmankac}
\forall \ij \mbox{ and }s\in [t,T], Y^{ij}_s=\bar v^{ij}(s,\xtxs).
\ee
In this Markovian framework of randomness, this result improves substancially the one by Hu-Tang \cite{Hutang-2008} on the same subject. Uniqueness is even new.

Finally at the end of the paper, we deal with issues related to the link of the solution of systems 
(\ref{mainsyst-vi}) and (\ref{maxminsyst-vi}) and the value function of the zero-sum switching game. In some particular cases, we show that they are equal and a saddle-point for the game is obtained. 

To the best of our knowledge, these issues have not been addressed in the literature yet. 
\ms

The paper is organized as follows. In Section 2, we fix some notations and, for sake of completeness, recall accurately under which conditions each of the systems 
(\ref{mainsyst-vi}) and (\ref{maxminsyst-vi}) has a unique solution. Note that these results are already given in \cite{Djehichehammorlais14}. In Section 3, we show that if mainly the switching costs $\bar g_{ij}$, $\ij$, are ${\cal C}^{1,2}$ then the unique solutions of (\ref{mainsyst-vi}) and (\ref{maxminsyst-vi}) coincide. Next the link with the unique solution of system (\ref{general-drbsde}) is  stated. The proof of this result is postponed to Appendix, given in the end of the paper. In Section 4, we first describe the zero-sum switching game problem. Then, under some additional conditions on the two families ($f^{ij}$)$_{\ij}$ and ($h^{ij}$)$_{\ij}$, we show that this game has a value which is given by the unique solution of (\ref{mainsyst-vi}) and (\ref{maxminsyst-vi}) and thus in terms of the solution of (\ref{general-drbsde}) as well, due to relation (\ref{Feynmankac}). We also provide a saddle-point of this game. The relationship of $Y^{ij}$ and the upper and lower values of the game are also studied. 

\section{Notations and first results}
Let $T$ (resp. $k$, $d$) be a fixed positive constant (resp.
two integers) and $\Gamma^{1}$ (resp. $\Gamma^{2} $) denote the set
of switching modes for player 1 (resp. 2). For later use, we shall
denote by $\Lambda$ the cardinal of the product set
$\Gamma:=\Gamma^{1}\times \Gamma^{2} $ and for $(i,j)\in \G$,
$\gami:=\Gamma^1-\{i\}$ and $\gamj:=\Gamma^2-\{j\}$. For
$\vec{y}=(y^{kl})_{(k,l)\in \gam} \in \R^{\L}$, $\ij$, and
$\underline{y}\in \R$, we denote by $[(y^{kl})_{(k,l)\in \gam-\{i,j\}},\underline{y}]$
the matrix obtained from the matrix $\vec{y}=(y^{kl})_{\kl}$ by replacing the element
$y^{ij}$ with $\underline{y}$.

\ms\no For any $\ij$, let
$$
\begin{array}{l}
b:\, (t,x)\in [0,T]\times \R^k\mapsto b(t,x)\in \R^{k};\\
\sigma: \, (t,x)\in [0,T]\times \R^k\mapsto \sigma(t,x)\in \R^{k\times
d};\\ f^{ij}: \, (t,x,\vec{y},z)\in [0,T]\times
\R^{k+\Lambda+d}\mapsto
f^{ij}(t,x,\vec{y},z)\,\,\in \R\,;\\
\underline{g}_{ik}: \, (t,x)\in [0,T]\times\R^k\mapsto \underline{g}_{ik}(t,x)\in \R\,\, ;\\
\bar{g}_{jl}: \, (t,x)\in [0,T]\times\R^k \mapsto \bar{g}_{jl}(t,x)\in \R\,\,;\\
h^{ij}: \, x\in \R^k\mapsto  h^{ij}(x)\in \R.
\end{array}
$$
\noindent A function $\Phi: \, (t,x)\in \espo\mapsto
\Phi(t,x)\in \R$ is called of {\it polynomial growth} if there exist
two non-negative real constants $C$ and $\g$ such that
$$ |\Phi(t,x)|\leq
C(1+|x|^\g),\, \forall (t,x)\in \espo.
$$
Hereafter, this class of functions is denoted by
$\Pi_{g}$. Let $\cC^{1,2}(\espo)$ (or simply
$\cC^{1,2}$) denote the set of real-valued functions defined on $\espo$,
which are once (resp. twice) differentiable w.r.t. $t$ (resp. $x$)
and with continuous derivatives.

\ms The following assumptions (H0)-(H4) on the data of the systems (\ref{mainsyst-vi}) and (\ref{maxminsyst-vi})  are in force throughout the paper. They are the same as
 in \cite{Djehichehammorlais14}.
\begin{enumerate}
\item[$\mathbf{(H0)}$] The functions $b$ and $\sigma$ 
are jointly continuous in $(t,x)$ and Lipschitz continuous w.r.t. $x$ uniformly in $t$, meaning that there exists a non-negative constant $C$ such that for any $(t,x,x')
 \in  [0,T] \times \mathbb{R}^{k+k}$ we have
$$|\sigma(t,x)- \sigma(t,x')|+ |b(t,x)- b(t,x')|
\le C|x -x'|.
$$
Thus they are also of linear growth w.r.t. $x$, i.e., there exists a constant $C$ such that for any $\tx$, 
$$
 |b(t,x) |+|\sigma(t,x)| \le C(1 +|x|).
$$

\item[$\mathbf{(H1)}$] Each function $f^{ij}$

(i) is continuous in $(t,x)$ uniformly w.r.t. the other variables
$(\vec{y},z)$ and, for any $(t, x)$, the mapping $(t,x)
\rightarrow f^{ij}(t,x,0,0)$ is of polynomial growth.

(ii) is Lipschitz continuous with
respect to the variables ($\vec{y}:=(y^{ij})_{(i,j)\in
\Gamma_1\times \Gamma_2},z$) uniformly in $(t,x)$, i.e.
 $\forall \; (t,x) \in [0,T] \times \mathbb{R}^{k},\; \forall \; (\vec{y}_{1}, \vec{y}_{2}) \in \mathbb{R}^{\Lambda} \times  \mathbb{R}^{\Lambda}
, (z^{1}, z^{2}) \in \mathbb{R}^{d} \times \mathbb{R}^{d},$
$$ |f^{ij}(t,x, \vec{y}_{1},z_{1}) -f^{ij}(t,x, \vec{y}_{2},z_{2}) | \le C\left( |\vec{y}_{1}-\vec{y}_{2}| +|z_{1}-z_{2}|\right),$$
where $|\vec{y}|$ stands for the standard Euclidean norm of
$\vec{y}$ in $\R^\L$.
\item[$\mathbf{(H2)}$] \underline{Monotonicity}: Let $\vec{y} = (y^{kl})_{(k,l) \in \Gamma^{1}\times \Gamma^{2}} $. For any $\ij$ and any $(k,l) \neq (i,j)$ the mapping $y^{kl}
\rightarrow f^{ij}(s,\vec{y}, z)$ is non-decreasing.
\item[$\mathbf{(H3)}$] The functions $h^{ij}$, which are the terminal conditions in the systems (\ref{mainsyst-vi}) and (\ref{maxminsyst-vi}), are continuous with respect to $x$, belong to class $ \Pi_{g}$ and satisfy the following consistency condition:
$$ \displaystyle{ \forall \; (i,j) \in \Gamma^{1} \times \Gamma^{2} \mbox{ and }x\in \R^k, \,\, \max_{k\in
(\Gamma^{1})^{-i}}\big( h^{kj}(x) -\underline{g}_{ik}(T,x) \big)
\le h^{ij}(x)  \le \min_{l\in (\Gamma^{2})^{-j}} \big(h^{il}(x)
+\bar{g}_{jl}(T,x)\big).} $$
\item[$\mathbf{(H4)}$] \underline{The non free loop property}: The switching costs $\underline{g}_{ik} $ and
$ \bar{g}_{jl} $ are non-negative, jointly continuous in $(t,x)$,
belong to $\Pi_g$ and satisfy the following condition:

For any loop in $\gam$, i.e., any sequence of pairs $(i_1,j_1),\ldots,(i_N,j_N)$
 of $\gam$ such that $(i_N,j_N)=(i_1,j_1)$, \mbox{card}$\{
(i_1,j_1),\ldots,(i_N,j_N)\}=N-1$ and any $q=1,\ldots,N-1$, either
$i_{q+1}=i_q$ or $j_{q+1}=j_q$, we have: $\forall (t,x)\in \espo$,
\be\label{nonfreeloop3} \sum_{q=1,N-1}\varphi_{i_qj_{q}}(t,x)\neq
0, \ee where, $\forall \,\,\,q=1,\ldots,N-1,\,\,
\varphi_{i_qj_q}(t,x)=-\underline{g}_{i_qi_{q+1}}(t,x)\ind_{i_q\neq
i_{q+1}}+\bar{g}_{j_qj_{q+1}}(t,x)\ind_{j_q\neq j_{q+1}}$.
\end{enumerate}

This assumption implies in particular that \be \label{lp1}
\displaystyle{\forall \; (i_{1}, \ldots, i_{N})\in (\Gamma^1)^N
\;\textrm{such that } \;i_{N} =i_{1}\mbox{ and } \mbox{card}\{i_{1},
\ldots, i_{N}\}=N-1,\,\, \sum_{p=1}^{N-1}\underline{g}_{i_{k}
i_{k+1}}
> 0 } \ee and \be \label{lp2} \displaystyle{\forall \; (j_{1},
\ldots, j_{N})\in (\Gamma^2)^N \;\textrm{such that } j_{N} =j_{1}
\mbox{ and } \mbox{card}\{j_{1}, \ldots, j_{N}\}=N-1, \;
\sum_{p=1}^{N-1}\bar{g}_{j_{k} j_{k+1}} > 0}.\ee By convention we
set $\bar{g}_{jj} = \underline{g}_{ii} =0$.

\ms\no Conditions (\ref{lp1})
and (\ref{lp2}) are classical in the literature of switching
problems and usually referred to as the {\it non free loop property}.
\ms

We now introduce the probabilistic tools we need later. Let $(\Omega, {\cal F}, \P)$ be a fixed probability space on which
is defined a standard $d$-dimensional Brownian motion
$B=(B_t)_{0\leq t\leq T}$ whose natural filtration is\\
$(\cF_t^0:=\sigma \{B_s, s\leq t\})_{0\leq t\leq T}$. Let $
\bF=(\cF_t)_{0\leq t\leq T}$ be the completed filtration of
$(\cF_t^0)_{0\leq t\leq T}$ with the $\mathbb{P}$-null sets of
${\cal F}$, hence $(\cF_t)_{0\leq t\leq T}$ satisfies the usual
conditions, i.e., it is right continuous and complete. On the other hand, 
we denote by ${\cal P}$ be the $\sigma$-algebra on $[0,T]\times \Omega$ of
$\bF$-progressively measurable sets. 

\ms\no  Next, let us fix $t$ in $[0,T]$ and let us define

(i) ${\cal H}_t^{2,\ell}$ ($\ell \geq 1$) be the set of $\cal
P$-measurable and $\R^\ell$-valued processes $w=(w_s)_{ s\in [0,T]}$
such that  \\$\E[\int_t^T|w_s|^2ds]<\infty$;

(ii) ${\cal S}_t^{2}$ (resp. ${\cal S}^2_{t, d}$) be the set of ${\cal P}$-measurable continuous (resp. RCLL) processes 
$w=(w_s)_{ s\in [0,T]}$ such that \\ $\E[\sup_{ t \le s\le T}|{w}_s|^2]<\infty$ ;

(iii) $\mathcal{A}^2_{t, i}$ be the subset
of ${\cal S}_t^{2}$ of non-decreasing processes $K=(K_s)_{ s\in [0, T]}$
such that $K_t=0$ (and then $K_s=0$ for $s\leq t$) ;

(iv) The sets ${\cal H}_0^{2,\ell}$,  ${\cal S}_0^{2}$, ${\cal S}^2_{0, d}$ and $\mathcal{A}^2_{0, i}$ will be simply denoted by ${\cal H}^{2,\ell}$,  ${\cal S}^{2}$, ${\cal S}^2_{d}$ and $\mathcal{A}^2_{ i}$.
\medskip

\no For $\tx$, let  $X^{t,x}$ be the diffusion process solution of the following standard SDE:
\be\label{sdex}\begin{array}{l}
\forall s\in [t,T],\,\,\, X_s^{t,x}=x+\int_t^sb(r,X^{t,x}_r)dr+\int_t^s\si(r,X^{t,x}_r)dB_r;\,\,\,X_s^{t,x}=x \mbox{ for }s\in [0,t].\end{array}
\ee
Under Assumption (H0) on $b$ and $\si$, the process $X^{t,x}$ exists and is unique (\cite{revuzyor}, Theorem 2.1 pp.375). Moreover, it satisfies the following estimates: For all $p\geq 1$,
\be\label{estimationx}
\E[\sup_{s\leq T}|X_s^{t,x}|^p]\leq C(1+|x|^p).
\ee
Its infinitesimal generator $\cL^X$ is given, for every $\tx$ and $\phi \in {\cC}^{1,2}$, by
\be\label{generateur}\begin{array}{lll}\cL^X\phi(t,x)&:=&\frac{1}{2}\sum\limits_{i,j=1}^k(\sigma\sigma^\top (t,x))_{i,j}\partial^2_{x_ix_j}\phi(t,x)+\sum_{i=1,k}b_i(t,x)\partial_{x_i}\phi(t,x)\\{}&\,\,=&
\frac{1}{2}Tr[\si\si^\top(t,x)D_{xx}^2\phi(t,x)]+b(t,x)^\top D_x\phi(t,x).\end{array}\ee
Under Assumptions (H0)-(H4), we have

\begin{theorem}\label{thmdhm} {(\cite{Djehichehammorlais14}, Theorems 5.4 and 5.5)} There exists a unique continuous viscosity solution in the class $\pl$ $(\bar v^{ij})_{\ij}$ (resp. 
$(\udl v^{ij})_{\ij}$) of the following system: $\forall \ij$,
\begin{eqnarray}\label{system1} 
         \left\{ \begin{array}{l} 
   \min\Big\{(\bar v^{ij }-L^{ij}(\vec {\bar  v}))(t,x);\max\Big\{ (\bar v^{ij} -U^{ij}(\vec {\bar v}))(t,x); \\\qq\qq\qq-\partial_t \bar v^{ij} (t,x)-\mathcal{L}^{X}(\bar v^{ij})(t,x) -f^{ij}(t,x, (\bar v^{kl}(t,x))_{\kl},\sigma^\top(t,x)D_x\bar v^{ij}(t,x)) \Big\}  \Big\}  =0,\\
   \bar v^{ij}(T,x)=h^{ij}(x)\\
           \end{array}          \right. \end{eqnarray}
(resp.
\begin{eqnarray}\label{system2}
         \left\{ \begin{array}{l} 
   \max\Big\{(\udl v^{ij} -U^{ij}(\vec{\udl  v}))(t,x); \min\Big\{(\udl  v^{ij }(t,x)-L^{ij}(\vec{\udl v}))(t,x);\\\qq   -\partial_t \udl v^{ij}(t,x) -\mathcal{L}^{X}(\udl v^{ij})(t,x) -f^{ij}(t,x, (\udl v^{kl}(t,x))_{\kl},\sigma(t,x)^\top D_x \udl v^{ij}(t,x)) \Big\}  \Big\}  =0,\\
  \udl   v^{ij}(T,x)=h^{ij}(x))  \\
           \end{array}          \right. \end{eqnarray}
where the obstacles $U^{ij}$ and $L^{ij}$ are defined in (\ref{expbar}). $\qed$
\end{theorem}
In order to obtain the solutions of the systems (\ref{system1}) and (\ref{system2}) respectively, Djehiche et al.
(\cite{Djehichehammorlais14}) introduced the following sequences of backward reflected BSDEs with inter-connected obstacles: $\forall m,n\geq 0$, $\forall
(i,j)\in \gam$, \be\label{penalizedscheme}
\left\{\begin{array}{l}
\bar Y^{ij,m}\in \cS^{2}, \,\,\bar Z^{ij,m} \in \cH^{2,d} \mbox{
and
}\bar K^{ij,m } \in \cA^{2}_{i}\,\,;\\
\bar Y^{ij,m}_s=h^{ij}(X^{t,x}_T)+\int_{s}^{T}\bar
{f}^{ij,m}(r,X^{t,x}_r,
(\bar Y^{kl,m}_r)_{(k,l)\in \gam},\bar Z^{ij,m}_r)dr+\int_s^Td\bar K^{ij,m}_r-\int_s^T\bar Z^{ij,m}_rdB_r,\,s\leq T;\\
\displaystyle{\bar Y^{ij,m}_{s}\geq \max_{k\in
(\Gamma^{1})^{-i}}\{\bar
Y^{kj,m}_s-\underline{g}_{ik}(s,X^{t,x}_s)\}},\,s\leq T;
\\
\int_0^T(\bar Y^{ij,m}_s-\max_{k\in (\Gamma^{1})^{-i}}\{\bar
Y^{kj,m}_s-\underline{g}_{ik}(s,X^{t,x}_s)\})d\bar
K^{ij,m}_s=0\end{array}\right.\ee and
\begin{equation} \label{second-penalizedscheme}  \left\{ \begin{array}{l} 
\underbar Y^{ij,n}\in \cS^{2}, \,\,\underbar Z^{ij,n} \in \cH^{2,d} \mbox{
and }\underbar K^{ij,n } \in \cA^{2}_{i}\,\,;\\
    \underbar Y_{s}^{ij,n} = h^{ij}(X_{T}^{t,x})+\int_{s}^{T}\underbar f^{ij,n}(r, X_r^{t,x}, (\underbar Y^{kl,n}_r)_{(k,l)\in \gam},\underbar Z^{ij,n}_r)dr - \int_s^T\underbar Z_r^{ij,n}dB_r -\int_{s}^{T}d\underbar K_r^{ij,n},   \,s\leq T;\\
\underbar Y^{ij,n }_s \le \min_{l\in \gamj} \big(\underbar Y_s^{il,n} +\bar{g}_{jl}(s,X_s^{t,x})\big),\,s\leq T; \\
\int_0^T(\underbar Y^{ij,n }_s - \min_{l\in \gamj} \{\underbar Y_s^{il,n} +\bar{g}_{jl}(s,X_s^{t,x})\})d\underbar K^{ij,n}_s =0 \end{array}\right.
\end{equation}
where, for any $(i,j) \in \gam$, $n,m\geq 0$ and $(s, x,\vec{y}, z^{ij})$,

\begin{equation}\label{f-max}\bar{f}^{ij,m}(s, x,\vec{y},z^{ij}):= f^{ij}(s,x,(y^{kl})_{(k,l)\in \gam},z^{ij})
-m\big(y^{ij} -\min_{l\in \gamj} (y^{il} + \bar {g}_{jl}
(s,x))\big)^{+}
\end{equation}
and 
\begin{equation}\label{f-min}
\underbar {f}^{ij,n}(s, x,\vec{y},z^{ij}):= f^{ij}(s,x,(y^{kl})_{(k,l)\in \gam},z^{ij})
+n\big(y^{ij} -\max_{k\in \gami} (y^{kj} - \underline{g}_{ik}
(s,x))\big)^{-}.
\end{equation}
Under Assumptions (H0)-(H4), it is shown in \cite{Hammorlais13} (see also \cite{cek2} or \cite{Hamzhang}) that each one of the systems (\ref{penalizedscheme}) and (\ref{second-penalizedscheme}) has a unique solution $(\bar Y^{ij,m},\bar Z^{ij,m},\bar K^{ij,m})$ and $(\underbar Y^{ij,m},\underbar Z^{ij,m},\underbar K^{ij,m})$ respectively. In addition, they enjoy the following properties:
\ms

\no (i) For any $m,n\geq 0$ and $(i,j)\in \gam$
\be \label{ineqynm}
\bar Y^{ij,m}\geq \bar Y^{ij,m+1}\geq \underbar Y^{ij,n+1}\geq \underbar Y^{ij,n}.
\ee
(ii) For any $n,m\geq 0$ and $(i,j)\in \gam$, there exist deterministic continuous functions $\bar v^{ij,m}$ and $\udl v^{ij,n}$ such that, for any $\tx$ and $s\in [t,T]$, we have
$$\bar Y^{ij,m}_s=\bar v^{ij,m}(s,X^{t,x}_s)\,\,\, \mbox{ and } \,\,\, \underbar Y^{ij,n}_s=\udl v^{ij,n}(s,X^{t,x}_s).$$
Moreover, from (\ref{ineqynm}) we easily deduce that, for any $n,m\geq 0$ and $\ij$,
\be\label{inegvijnm}
\bar v^{ij,m}\ge \bar v^{ij,m+1}\ge  \udl v^{ij,n+1}\ge \udl v^{ij,n}.
\ee
Finally, for any $m\ge 0$ (resp. $n\ge 0$), $\bar  v_m:=(\bar v^{ij,m})_{\ij}$ (resp.
$\udl v_n:=(\udl v^{ij,n})_{\ij}$) is the unique continuous viscosity solution, in the class $\pl$, of the following system of PDEs with inter-connected obstacles: $\forall \ij$, $\forall \tx$,
\begin{eqnarray*} 
         \left\{ \begin{array}{l} 
   \min\Big\{(\bar v^{ij,m }-L^{ij}(\vec {\bar v}_m)(t,x);\\  \qq -\partial_t \bar v^{ij,m}(t,x) -\mathcal{L}^{X}(\bar v^{ij,m})(t,x) -\bar f^{ij,m}(t,x, (\bar v^{kl,m}(t,x))_{\kl},\sigma(t,x)^\top D_x \bar v^{ij,m}(t,x)) \Big\} =0,\\
   \bar  v^{ij,m}(T,x)=h^{ij}(x) \\
           \end{array}          \right. \end{eqnarray*}
(resp.
\begin{eqnarray*} 
         \left\{ \begin{array}{l} 
   \max\Big\{(\udl v^{ij,n} -U^{ij}(\vec{\udl v}_n))(t,x); \\\qq-\partial_t \udl v^{ij,n}(t,x) -\mathcal{L}^{X}(\udl v^{ij,n})(t,x) -\underbar f^{ij,n}(t,x, (\udl v^{kl,n}(t,x))_{\kl},\sigma^\top(t,x)D_x\udl v^{ij,n}(t,x)) \Big\}  =0,\\
   \udl  v^{ij,n}(T,x)=h^{ij}(x)). \\
           \end{array}          \right. \end{eqnarray*}
(iii) For $\ij$ and $(t,x)\in \espo$, let us set 
$$
\bar v^{ij}(t,x):=\lim_{m\rw \infty}\searrow \bar v^{ij,m}(t,x)\;\;\; \mbox{ and } \;\;\; 
\udl v^{ij}(t,x):=\lim_{n\rw \infty}
\nearrow \udl v^{ij,n}(t,x).
$$ 
Then, using Perron's method, it is shown that $(\bar v^{ij})_{\ij}$ (resp. 
$(\udl v^{ij})_{\ij}$) is continuous, belongs to $\pl$ and is the unique viscosity solution, in class $\pl$, of system (\ref{system1}) (resp. (\ref{system2})). Finally, by construction and in view of (\ref{inegvijnm}), it holds that, for any $\ij$,
\be\label{ineqvij}\udl v^{ij}\le \bar v^{ij}.\ee      

\section{Equality of the solutions of min-max and max-min systems. Related system of reflected BSDEs.}
In \cite{Djehichehammorlais14}, the question whether or not  for any $\ij$, $\udl v^{ij}\equiv \bar v^{ij}$ was left open. This was mainly due to the fact we have not been able to compare the inter-connected obstacles neither in (\ref{system1}) nor in (\ref{system2}). 
\ms

Actually, had we known that 
\be\label{inegalite}
\begin{array}{l}
\mbox{ (i) }\q \forall \ij, \;\;\; L^{ij}(\vec{\bar v})\leq U^{ij}(\vec{\bar v})\\ \mbox{or}\\\mbox{ (ii) }\q\forall \ij, \;\;\; L^{ij}(\vec{\udl v})\leq U^{ij}(\vec{\udl v})\ea
\ee
then we would have deduced, from the general existence result obtained in 
Hamad\`ene-Hassani \cite{Hamhassani}  (see Theorem \ref{existence} in appendix (A2)) and the uniqueness of the solution of (\ref{system1}) or 
(\ref{system2}), that for any $\ij$, $\bar v^{ij}=\udl v^{ij}$.
In this section, we are going to investigate under which additional regularity assumptions on the data of the problem, one of the inequalities in (\ref{inegalite}) is satisfied to be able to conclude that $\bar v^{ij}=\udl v^{ij}$, for any $\ij$, i.e., the solutions of (\ref{system1}) and (\ref{system2}) are the same.

\ms \no For this objective, let us introduce the following additional assumption.

\noindent {\bf (H5)}: 
\ms

\no (i) For any $(i,j)\in \gam$, the functions $\bar g_{ij}$ are ${\cal C}^{1,2}$ and, $D_x\bar g_{ij}$, $D^2_{xx}\bar g_{ij}$ belong to $\pl$. Furthermore, for any $j_1,j_2,j_3 \in \Gamma_2$ such that $|\{j_1,j_2,j_3\}|=3$,  
$$
\bar g_{j_1j_3}(s,x)<\bar g_{j_1j_2}(s,x)+\bar g_{j_2j_3}(s,x),\,\,\forall (s,x)\in [0,T]\times \R^k.
$$ 
(ii) For any $\ij$, the function $f^{ij}$ verifies the following estimate:
$$
|f^{ij}(s,x,\vec{y},z^{ij})|\leq C(1+|x|^p),\,\,\forall (s,x,\vec{y},z^{ij})\in [0,T]\times \R^{k+\Lambda+d},$$for some real constants $C$ and $p$.
\begin{rem}\label{penalities} Note that by It\^o's formula, for any $\ij$, 
$$ \left\{\begin{array}{l}
\bar g_{ij}(s,X^{t,x}_s)=\bar g_{ij}(t,x)+\int_t^s{\cal L}^X(\bar g_{ij})(r,X^{t,x}_r)dr+\int_t^sD_x\bar g_{ij}(r,X^{t,x}_r)\sigma(r,X^{t,x}_r)dB_r, \,\;\; s\in [t,T]\\ 
\mbox{and}\\
\bar  g_{ij}(s,X^{t,x}_s)=\bar g_{ij}(s,x),\;\;\; s\leq t.\end{array}
 \right.$$
Hereafter and to ease the reading
of the It\^o-Tanaka formula in Step 2 below (proof of Theorem \ref{thmprincipalsection1}), we denote by $\alpha^{ij}$ and $\beta^{ij}$, $\ij$, the following processes:
$$
\alpha^{ij}(s):={\cal L}^X(\bar g_{ij})(s,X^{t,x}_s),\,\, \beta^{ij}(s):= D_x\bar g_{ij}(s,X^{t,x}_s)\sigma(s,X^{t,x}_s),\,\, s\leq T. 
$$ 
\end{rem}

\ms We now provide the main result of this section.  
\begin{thm}\label{thmprincipalsection1}
Under Assumptions (H0)-(H5), for any $\ij$, it holds that
$$\bar v^{ij}=\udl v^{ij}.\qq\qq\qq\qq\qq\qq $$
\end{thm}
 \noindent We derive this last equality after the following four steps.
\ms

\no \udl{{\bf Step 1}}: Another approximating scheme for system (\ref{system1}).
\ms

For any $m\geq 0$, $\ij$ and $\tx$, let us consider the system of reflected BSDEs with one interconnected obstacle:
\be\label{intermpenscheme}
\left\{\begin{array}{l}
Y^{ij,m}\in \cS^{2}, \,\, Z^{ij,m} \in \cH^{2,d} \mbox{
and
} K^{ij,m } \in \cA_{i}^{2}\,\,;\\
 Y^{ij,m}_s=h^{ij}(X^{t,x}_T)+\int_{s}^{T}
{f}^{ij,m}(r,X^{t,x}_r,
( Y^{kl,m}_r)_{(k,l)\in \gam}, Z^{ij,m}_r)dr+\int_s^Td K^{ij,m}_r-\int_s^T Z^{ij,m}_rdB_r,\,s\leq T;\\
\displaystyle{ Y^{ij,m}_{s}\geq \max_{k\in
(\Gamma^{1})^{-i}}\{
Y^{kj,m}_s-\underline{g}_{ik}(s,X^{t,x}_s)\}},
s\leq T;\\
\int_0^T( Y^{ij,m}_s-\max_{k\in (\Gamma^{1})^{-i}}\{
Y^{kj,m}_s-\underline{g}_{ik}(s,X^{t,x}_s)\})d
K^{ij,m}_s=0,\end{array}\right.\ee 
where, 
\be \label{f-sum}
{f}^{ij,m}(s,x,\vec y,z^{ij}):={f}^{ij}(s,x,\vec y,z^{ij})-m\sum_{l\in \gamj}(y^{ij}-y^{il}-\bar g_{jl}(s,x))^+.
\ee
This generator is (slightly) different from $ \bar f^{ij, m}$ given by (\ref{f-max}) in Section 2. We mention that this new penalized generator is more convenient both for the application of the It\^o-Tanaka formula and obtention of the estimate 
(\ref{eq:bound-penalizedterm}) (in Step 2). On the other hand, note that for any $\ij$ and $(k,l)\neq (i,j)$ the mapping $\tilde y\in \R \mapsto {f}^{ij,m}(s,x,[(y^{rq})_{(r,q)\in \gam-(k,l)},\tilde {y}],z^{ij})$ is non-decreasing. \\
By Corollary 2, in \cite{Hammorlais13}, the solution of this system exists and is unique and there exist deterministic continuous functions $( v^{ij,m})_{\ij}$, which belong also to $\pl$ such that, for any $i,j$ and $m\geq 0$, it holds that
\be \label{eqyijm}
\forall s\in [t,T],\quad  Y^{ij,m}_s=v^{ij,m}(s,X^{t,x}_s).
\ee
Moreover, the family of functions $ \vec v_m:=( v^{ij,m})_{\ij}$ is the unique continuous solution in viscosity sense in $\pl$ of the following system of PDEs with obstacles:
\begin{eqnarray*} 
         \left\{\begin{array}{l} 
   \min\Big\{( v^{ij, m }-L^{ij}(\vec {v}_m))(t,x) ; \\\qq\qq-\partial_t  v^{ij, m}(t,x) -\mathcal{L}^{X}( v^{ij, m})(t,x) -f^{ij,m}(t,x, ( v^{kl,m}(t,x))_{\kl},\sigma(t,x)^\top D_x v^{ij, m}(t,x)) \Big\} =0,\\
     v^{ij,m}(T,x)=h^{ij}(x).\\
\end{array}\right. \end{eqnarray*}
Finally, by the Comparison Theorem (see \cite{Hammorlais13}, Remark 1) and using that
$ f^{ij,m+1}\le  f^{ij,m}$ and  $\bar f^{ij, |\Gamma_2|m}\leq  f^{ij,m}\leq \bar f^{ij,m}$, we deduce:
 $\forall \ij$ and $m\geq 0$, 
$$
 Y^{ij,m+1}\leq Y^{ij,m}\quad \mbox{ and }\quad \bar Y^{ij,|\Gamma_2| m}\leq  Y^{ij,m}\leq \bar Y^{ij,m},
$$
which implies that, for any $\ij$ and $m\geq 0$,
$$
 v^{ij,m+1}\leq v^{ij,m} \mbox{ and }\bar v^{ij,|\Gamma_2| m}\leq 
 v^{ij,m}\leq \bar v^{ij,m}.
$$
Then, for any $\ij$, the sequence $( v^{ij,m})_{m\geq 0}$ is decreasing and converges, by Dini's theorem, uniformly on compact subsets of $\espo$, to $\bar v^{ij}$ since $\lim_{m\rw \infty}\bar v^{ij,m}(t,x)=\bar v^{ij}(t,x)$, for any $(t,x)\in \espo$.
\bs

\no \udl{\bf Step 2}: The following estimate holds true: For any $t\leq T$, $\ij$ and $m\ge 0$, 
\begin{equation}\label{eq:bound-penalizedterm}\begin{array}{l}
\mathbb{E}\Big\{ m\int_{t}^{T}\sum_{l\in \gamj}\{ Y^{ij,m}_s - Y^{il,m}_{s} -\bar{g}_{jl}(s, X_s^{t,x})\}^+ds \Big \} \le C(1+|x|^p),\end{array}
\end{equation}
and 
\begin{equation}\label{eq:bound-penalizedterm2}\begin{array}{l}
\mathbb{E}\Big\{ m^2\int_{t}^{T}\sum_{l\in \gamj}(\{ Y^{ij,m}_s - Y^{il,m}_{s} -\bar{g}_{jl}(s, X_s^{t,x})\}^+)^2 ds \Big \} \le C(1+|x|^{2p}),\end{array}
\end{equation} 
where $p$ and the generic constant $C$ are independent of $m$ and $x$.
\bs

\no For later use, we first give a representation of $ Y^{ij,m}$ as the optimal payoff of a switching problem. Indeed,  let $\d:=(\t_n,\zeta_n)_{n\geq 0}$ be an admissible strategy of switching, i.e.,   

(a) $(\t_n)_{n\ge 0}$ is an increasing sequence of stopping times such that $\,\, \P[\t_n<T, \forall n\geq 0]=0$;

(b) $\forall n\ge 0$, $\zeta_n$ is a random variable with values in $\Gamma^1$ and $\cF_{\t_n}$-measurable;  

(c) If $(A^\d_s)_{s\leq T}$ is the non-decreasing, $\bF$-adapted and RCLL process defined by
$$
\forall s\in [0,T),\quad A^\d_s=\sum_{n\geq 1}\underline g_{\zeta_{n-1}\zeta_n}(\t_n,X^{t,x}_{\t_n})\ind_{\{\t_n\leq s\}}\quad  \mbox { and }\quad A_T^\d=\lim_{s\rw T}A^\d_s,
$$ then $\E[(A^\d_T)^2]<\infty$. The quantity $A^\d_T$ stands for the switching cost at terminal time $T$ when the strategy $\d$ is implemented.

\no Next, with an admissible strategy $\d:=(\t_n,\z_n)_{n\ge 0}$ we associate a piecewise constant process 
$a=(a_s)_{s\in [0,T]}$ defined by \be a_s:=\z_0\mathbf{1}_{\{\t_0\}}(s)+\sum\limits^\infty_{j=1}\z_{j-1}\mathbf{1}_{]\t_{j-1},\t_j]}(s),\quad  s\leq T.\ee
For any $s\geq \t_0$, $a_s$ is the mode indicator at time $s$ of the system which is subject to control strategy $\d$. Note that there is a bijection between the processes $a$ and the admissible strategies $\d$, therefore hereafter we indifferently  write $A^a$ or $A^\d$. This notation $a$ for the indicator process shall be used in Section 4 to deal with the zero-sum switching game (see Theorem \ref{gamevalue}).

\no Finally, for any fixed $i\in \Gamma^1$ and a real constant $\th \in [t,T]$, we denote by $\cA^i_\th$ the following set:
$$
\cA^i_\th:=\Big \{\d=(\t_n,\z_n)_{n\ge 0} \mbox{ admissible strategy such that } \t_0=\th \mbox{ and }\z_0=i\Big \}.
$$

\no Now, for any $\d=(\t_n,\a_n)_{n\ge 0}$ (or equivalently $a$) which belongs to $\cA^i_\theta$, let us define the pair of processes  
$( U^{aj,m}, V^{aj,m})$ which belongs to $\cS^2_{d}\times \cH^{2,d}$ and which
 solves the following BSDE (which is of non standard form):  
\begin{equation} \label{switchedbsdebis} 
  U_{s}^{aj,m} = h^{a_T j}(X_{T})  +\int_{s}^T 1_{\{r\geq \t_0\}} f^{aj,m}(r, X^{t,x}_r, U_r^{aj,m}, V^{aj,m}_{r})dr -\int_{s}^{T} V^{aj,m}_{r}dB_r - \big(A_{T}^{a} -A_{s}^{a}\big),\,\, s\le T,
\end{equation}
where, for any $s\geq \t_0$ and $(\bar y,\bar z)\in \R^{1+d}$, 
$ f^{aj,m}(s,X^{t,x}_s, \bar y,\bar z)$ (resp. $ f^{aj}(s,X^{t,x}_s,\bar y,\bar z)$) is equal to $$ f^{\ell j,m}(s,X^{t,x}_s,[( v^{kl,m}(s,X^{t,x}_s))_{(k,l)\in \gam -\{(\ell ,j)\}},\bar y],\bar z)$$(resp. $$ f^{\ell j}(s,X^{t,x}_s,[(v^{kl,m}(s,X^{t,x}_s))_{(k,l)\in \gam -\{(\ell ,j)\}},\bar y],\bar z))$$
if at time $s$, $a(s)=\ell$. 
Let us point out that since $a$ is admissible and then $\E[(A^\d_T)^2]<\infty$, the solution of equation (\ref{switchedbsdebis}) exists and is unique by an immediate change of variables.  Furthermore, we have the following representation of $ Y^{ij,m}$ (see e.g. \cite{Hammorlais13,Hutang} for more details on this representation):
\begin{equation} \label{switchingrepbis}
\dis{ Y_\th^{ij,m}=  \esssup_{a \in \mathcal{A}_{\th}^i}  \{ U_{\th}^{a,j,m} -A^a_\th\}},\, t \le \theta\leq T.
\end{equation}

The equality (\ref{switchingrepbis}) differs from the one given in \cite{Hutang} and some other papers including \cite{Hammorlais13}. However there is a lack in the previous papers which we correct here. Note that this is a minor point which does not affect the results in those papers (\cite{cek2, Hammorlais13,Hutang} etc.). The accurate relation is given in (\cite{Djehichehampopier}, equation (9), pp. 2757) in the particular case when the generators do not depend on the components $(\vec{y},z)$ but this fact is irrelevant. 

Finally note that the function $ f^{\ell j,m}$ depends only on $(\bar y,\bar z)$. However the representation (\ref{switchingrepbis}) for $Y^{ij,m}$ still holds since the solution of system of reflected BSDEs (\ref{intermpenscheme}) is unique and by (\ref{eqyijm}). It follows that, for any $j,l\in \Gamma^2$ and $\th \le T$,  
\begin{equation}\label{domination}( Y^{ij,m}_\th- Y_\th^{il,m} -\bar{g}_{jl}(\th,X^{t,x}_\th))^+\le \esssup_{a \in \mathcal{A}_{\th}^i}  (U_{\th}^{aj,m} - U_{\th}^{al,m} -\bar{g}_{jl}(\th,X^{t,x}_\th))^{+}.
\end{equation}
We now examine the quantity $( U_{\th}^{aj,m} - U_{\th}^{al,m} -\bar{g}_{jl}(\th,X^{t,x}_\th))^{+}$. Define the set 
$\mathcal{B}_{jl}$ as follows:
$$\mathcal{B}_{jl}= \{(s,\omega)\in [0,T]\times \Omega,\, \;\; \textrm{such that} \;\; U_s^{aj,m} - U_s^{al,m}-\bar{g}_{jl}(s,X^{t,x}_s)> 0\}$$ and, for any $s\in [0,T]$, 
\be \label{eqwijm}
W^{a,jl,m}_s:= U_s^{aj,m} - U_s^{al,m}-\bar{g}_{jl}(s,X^{t,x}_s).
\ee
\newpage
Then, by It\^o-Tanaka's formula, we have, for every $s\in [\th,T]$,
\[
\begin{array}{l}
(W^{a,jl,m}_{s})^{+} +\frac{1}{2}\int_{s}^{T} d{L}^{a,jl,m}_{r}+
m\int_{s}^{T}dr\{\sum_{j'' \neq j} 
 \mathbf{1}_{\mathcal{B}_{jl}(r)}(W^{a,jj^{\prime\prime},m}_r)^{+} 
-\sum_{j^{\prime\prime} \neq l } \mathbf{1}_{\mathcal{B}_{jl}}(r) (W^{a,lj^{\prime\prime},m}_r)^{+} \}
\\\\
\;\qq\qq=  \int_{s}^{T} \mathbf{1}_{\mathcal{B}_{jl}}(r) \{ f^{aj}(r,X^{t,x}_r, U_r^{aj,m}) - f^{al}(r,X^{t,x}_r, U_r^{al,m}) -\alpha^{jl}(r)\}dr \\\\
 \;\qq\qq \qq - \int_{s}^{T}\mathbf{1}_{\mathcal{B}_{jl}}(r)\big(V_r^{aj,m} - V_r^{al,m} - \beta^{jl}({r})\big) dB_{r}\\

\end{array}
\]
where, the process ${L}^{a,jl,m}$ is the local time at 0 of the semimartingale $W^{a,jl,m}$. Splitting the difference 
$$\dis{ \Delta_{a,jl,m }(r) := m \sum_{j^{\prime\prime} \neq j}  \mathbf{1}_{\mathcal{B}_{jl}}(r)
(W_{r}^{a,jj^{\prime\prime},m})^{+} -m \sum_{j^{\prime\prime} \neq l } \mathbf{1}_{\mathcal{B}_{jl}}(r)( W_{r}^{a,lj^{\prime\prime},m})^{+}}
$$
as  
$$ \dis{ \Delta_{a, jl,m }(r) = m\mathbf{1}_{\mathcal{B}_{jl}}(r)(W_{r}^{a,jl,m})^+
- \mathbf{1}_{\mathcal{B}_{jl}}(r) (W_{r}^{a,lj,m})^++ m \sum_{j^" \neq j,l }\mathbf{1}_{\mathcal{B}_{jl}}(r)\{(W_{r}^{a,jj^{\prime\prime},m})^+ - (W_{r}^{a,lj'',m})^+}\},$$
the previous formula can be rewritten as follows: $\forall s\in [\theta,T]$, 
\be\label{equationestim}
\begin{array}{l}
(W^{a,jl,m}_{s})^{+} +\frac{1}{2}\int_{s}^{T} d{L}^{a,jl,m}_{r}
+ m\int_s^T \mathbf{1}_{\mathcal{B}_{jl}}(r)(W_{r}^{a,jl,m})^+ dr
\\\\
\qq\qq=  \int_{s}^{T} \mathbf{1}_{\mathcal{B}_{jl}}(r) ( f^{aj}(r,X^{t,x}_r, U_r^{aj,m}, V_r^{aj,m}) - f^{al}(r,X^{t,x}_r, U_r^{al,m}, V_r^{al,m}) -\alpha^{jl}({r}))dr+m\int_s^T\mathbf{1}_{\mathcal{B}_{jl}}(r) (W_{r}^{a,lj,m})^+dr \\\\
 \;\qq\qq - \int_{s}^{T}\mathbf{1}_{\mathcal{B}_{jl}}(r)\big( V_r^{aj,m} - V_r^{al,m} - \beta^{jl}({r})\big) dB_{r}-
m\int_{s}^{T}dr\{\sum_{j'' \neq j,l} 
 \mathbf{1}_{\mathcal{B}_{jl}(r)}[(W^{a,j^{\prime\prime},m}_r)^{+} -
(W^{a,lj",m}_r)^{+} ]\}.
 \end{array}
\ee
But by (H5)-(i), one has $\bar{g}_{jl} (t,x)+ \bar{g}_{lj} (t,x)> \bar{g}_{jj} (t,x)= 0 $.  Thus, we obtain that, for every  $(t,x)\in \espo$,   
$$ \{y \in \R^m,\; y_{j} -y_{l} -\bar{g}_{jl}(t,x) \ge 0 \}\cap \{y \in \R^m, \; y_{l} -y_{j} -\bar{g}_{lj}(t,x) \ge 0\} = \emptyset, $$
from which we deduce that  
\be\label{majoration0}\mathbf{1}_{\mathcal{B}_{jl}}(r) (W_{r}^{a,lj,m})^+ = 0,\,\,\forall r\in [\th, T]. \ee
Relying next on the elementary inequality $ a^{+} -b^{+} \le (a-b)^{+}$ ($a,b\in \R$), it holds
\begin{equation}\label{majoration}
\mathbf{1}_{\mathcal{B}_{jl}(r)}[(W^{a,jj",m}_r)^{+} -
(W^{a,lj",m}_r)^{+}] \le  
\mathbf{1}_{\mathcal{B}_{jl}}(r)\left( U^{al,m}_r - U^{aj,m}_r -\bar{g}_{lj^{\prime\prime}}(r,X^{t,x}_r) + \bar{g}_{jj^{\prime\prime}}(r,X^{t,x}_r)\right)^{+}.
\end{equation}
Using here that the family of penalty costs satisfies: $\bar{g}_{jj^{"}} <\bar{g}_{jl}+ \bar{g}_{lj^{"}}  $ 
 we deduce that
$$\{y \in \mathbb{R}^{m},\; y_{j} -y_{l} -\bar{g}_{jl}(t,x) \ge 0\} \cap \{y \in \mathbb{R}^{m},\; y_{l} -y_{j} - \bar{g}_{lj^{\prime\prime}}(t,x)+\bar{g}_{jj^{\prime\prime}}(t,x)  \ge 0\}=
 \emptyset $$
which therefore yields 
\be\label{majoration1}\forall r\in [\th,T],\,\, \mathbf{1}_{\mathcal{B}_{jl}}(r)\left( U^{al,m}_r - U^{aj,m}_r -\bar{g}_{lj^{\prime\prime}}(r,X^{t,x}_r)+ \bar{g}_{jj^{\prime\prime}}(r,X^{t,x}_r)\right)^{+} = 0.\ee
Going back now to (\ref{equationestim}), applying It\^o's formula to $e^{-m{s}}(W^{a,jl,m}_{s})^{+}$ and taking into account of (\ref{majoration0}), (\ref{majoration}) and (\ref{majoration1}) to obtain: 
$\forall s\in [\th,T]$, 
\[
\begin{array}{l}
(W^{a,jl,m}_{s})^{+} \leq
\int_{s}^{T} \mathbf{1}_{\mathcal{B}_{jl}}(r) e^{-m(r-s)}( f^{aj}(r,X^{t,x}_r,  U_r^{aj,m}, V_r^{aj,m}) - f^{al}(r,X^{t,x}_r, U_r^{al,m}, V_r^{al,m}) -\alpha^{jl}({r}))dr\\\\
 \;\qq\qq - \int_{s}^{T}\mathbf{1}_{\mathcal{B}_{jl}}(r)e^{-m(r-s)}\big( V_r^{aj,m} -  V_r^{al,m} - \beta^{jl}({r})\big) dB_{r}.
 \end{array}
\]
 Making now use of the estimates given in Assumptions (H0)-(H5) (namely the polynomial growth of both the drivers 
 ${f}^{ij}$ and of the penalty costs $\bar{g}_{ij}$) and
 taking the conditional expectation, we obtain:
$\forall s\in [\th,T]$, 
\[
\begin{array}{ll}
(W^{a,jl,m}_{s})^{+}& \leq
C\E[\int_{s}^{T} \mathbf{1}_{\mathcal{B}_{jl}}(r) e^{-m(r-s)}(1+|X^{t,x}_r|^p)dr|\cF_s]\\\\
& \leq
\frac{C}{m}\E[(1+\sup_{r\leq T}|X^{t,x}_r|^p)|\cF_s].
 \end{array}
\]
Recall now (\ref{domination}) and (\ref{eqwijm}) to obtain
\be \label{estimess}
m( Y^{ij,m}_\th- Y_\th^{il,m} -\bar{g}_{jl}(\th,X^{t,x}_\th))^+\leq 
C\E[(1+\sup_{r\leq T}|X^{t,x}_r|^p)|\cF_\th].
\ee
Taking expectation in both hand-sides and integrating in $\theta \in [t,T]$ to obtain (\ref{eq:bound-penalizedterm}). Next by squaring each side of the previous inequality, taking expectation and finally using Doob's inequality (\cite{revuzyor}, pp.54) we obtain:\be\label{estimalpha2}
m^2\E\Big \{\sum_{l\neq j}(( Y^{ij,m}_\th- Y_\th^{il,m} -\bar{g}_{jl}(\th,X^{t,x}_\th))^+)^2\Big \}\leq 
C(1+|x|^{2p}),\quad \forall \th\leq T
\ee 
since  $X^{t,x}$ has moments of any order by (\ref{estimationx}). Now as $\th$ is arbitrary in $[t,T]$ then, once more by integration with respect to $\th$ in the previous inequality, we obtain (\ref{eq:bound-penalizedterm2}). 
\ms

\no \udl{{\bf Step 3}}: We now prove that for any 
$(t_0,x_0)\in [0,T] \times \mathbb{R}^k$ and $\ij$, 
\be \label{inegalites-barrieres}L^{ij}(\vec {\bar v})(t_{0}, x_{0}) \le \bar v^{ij}(t_{0}, x_{0}) \le U^{ij}(\vec {\bar v})(t_{0}, x_{0}).\ee
We just need to check the property for $t_0<T$ since by the consistency condition (see (H3)), those inequalities hold true for $t_0=T$.   
\ms 

\no We first claim that $\bar v^{ij}(t_0,x_0) \ge L^{ij}(\vec {\bar v})(t_0,x_0)$ holds. Indeed, by construction of the sequence $\vec{\bar v}_m:=(\bar v^{ij,m})_{\ij}$ (proof of Theorem \ref{thmdhm}, Section 2), one has $\bar v^{ij,m} (t_0,x_0)\ge L^{ij}(\vec{\bar  v}_m)(t_0,x_0)$. Therefore, taking the limit w.r.t. $m$, we obtain $\bar v^{ij}(t_0,x_0) \ge L^{ij}(\vec {\bar v})(t_0,x_0)$. 
\ms

\no We now show that $\bar v^{ij}(t_{0}, x_{0}) \le U^{ij}(\vec{\bar v})(t_{0}, x_{0})$. First, assume that $\bar v^{ij}(t_0, x_0) > L^{ij}(\vec{\bar v})(t_0, x_0)$. Then, relying on the viscosity subsolution property of $\bar v^{ij}$ yields 
$$\begin{array}{l}\min\Big\{(\bar v^{ij }-L^{ij}(\vec{\bar v}))(t_0, x_0); \max\Big\{ (\bar v^{ij} -U^{ij}(\vec{\bar v}))(t_0, x_0); \\\qq\qq-\partial_t \bar v^{ij} (t_0, x_0)-\mathcal{L}^{X}(\bar v^{ij})(t_0, x_0) -f^{ij}(t_0,x_0, (\bar v^{kl}(t_0, x_0))_{\kl},\sigma(t_0,x_0)^\top D_x \bar v^{ij}(t_0,x_0)) \Big\}  \Big \}  
\le  0,\end{array}
 $$
which implies that $$
\max\Big\{ (\bar v^{ij} -U^{ij}(\vec{\bar v}))(t_0,x_0); -\partial_t \bar v^{ij}(t_0,x_0) -\mathcal{L}^{X}(\bar v^{ij})(t_0,x_0) -f^{ij}(t_0,x_0, (\bar v^{kl}(t_0,x_0))_{(k,l)\in \gam}))\Big \}\le 0. 
$$
Hence, $(\bar v^{ij} -U^{ij}(\vec{\bar v}))(t_0,x_0)\le 0$.\\
\ms

\no Suppose now that, at $(t_0, x_0)$, we have: $\bar v^{ij}(t_0, x_0) = L^{ij}(\vec{\bar v})(t_0, x_0)$. Proceeding by contradiction we suppose in addition that
\begin{equation} \label{contradiction} \dis{\exists \;
\epsilon>0,\quad (\bar v^{ij} -U^{ij}(\vec{\bar v}))(t_0,x_0) > \epsilon}.\end{equation}
Using both the continuity of $ (t,x) \mapsto \bar v^{ij}(t,x)
$ and $ (t,x) \mapsto U^{ij}(\vec{\bar v})(t,x)$ as well as the uniform convergence on compact subsets of $(v^{ij,m})_{m\ge 0}$ to $\bar v^{ij}$, we claim that for some strictly positive $\rho $ and for $m_{0} $ large enough it holds that
$$ \forall m\ge m_{0},\,\, \forall (t,x) \in \mathcal{B}((t_0, x_0), \rho),\,\quad ( v^{ij,m}-U^{ij}(\vec{ v}_m))(t,x) \ge \frac{\epsilon}{2}, $$
with $\mathcal{B}((t_0, x_0), \rho) = \{(t,x)\in \espo\,\, \textrm{s.t.}\,\, \; |t-t_{0}|\le \rho,\; |x-x_{0}| \le \rho \}$.\\
Without loss of generality we can now assume $[t_{0}, t_{0}+\rho] \subset [t_{0}, T]$. Let $(t,x)\in \mathcal{B}((t_0, x_0), \rho)$. By the definition of $U^{ij}(\vec{ v}_m)$ and as $\Gamma^2$ is finite, there exists one index $l_0 \neq j$ (which may depend on $(t,x)$) such that $$ v^{ij,m}(t,x)-( v^{il_0,m}(t,x) +\bar{g}_{jl_0}(t,x)) \ge \frac{\epsilon}{2}.$$By summing over $l\in \gami$, we deduce that for any $(t,x)\in \mathcal{B}((t_0, x_0), \rho)$, 
\begin{equation} \label{newcontradict}
 \sum_{l\in \gamj} \left( v^{ij,m}- v^{il,m}-\bar{g}_{jl} \right)^{+} (t,x)\ge\{ v^{ij,m}-( v^{il_0,m} +\bar{g}_{jl_0})\} ^+(t,x)\ge \frac{\epsilon}{2}.
 \end{equation}
Let us now introduce the following stopping time $\tau_{X}$:
$$ 
\tau_{X} = \inf\{ s \ge t_0,\; X_{s}^{t_{0}, x_{0}} \not \in\mathcal{B}((t_0, x_0), \rho) \} \wedge \big(t_{0}+\rho\big).
$$
 We  then have, for all $m \ge m_0$,
\begin{equation} \label{minoration} 
 m\mathbb{E}\Big \{\int_{t_0}^{\tau_X} \sum_{l\neq j} \{ v^{ij,m}(s,X^{t_0,x_0}_{s}) -\big(
v^{il,m}(s,X^{t_0,x_0}_{s}) +\bar{g}_{jl}(s, X^{t_0,x_0}_s))\}^{+}ds \Big \}
 \ge m\frac{\epsilon}{2}
\mathbb{E}(\tau_{X}-t_{0})\to\infty,
 \end{equation}
 as $m\to\infty$. But, this is contradictory to (\ref{eq:bound-penalizedterm}). Then $\bar v^{ij}(t_0,x_0)\leq U^{ij}(\vec{\bar v}))(t_0,x_0)$ and the proof of the claim is complete. 
\ms

\no \udl{\bf Step 4}: Finally, using inequality (\ref{inegalites-barrieres}) and Theorem 6.2 in \cite{Hamhassani} 
 (Theorem \ref{existence} in appendix (A2)), we deduce that for any fixed $\ij$, $\bar v^{ij}$ is continuous and of polynomial growth and also a viscosity solution of 
\begin{eqnarray}\label{system2.1}
         \left\{ \begin{array}{l} 
   \max\Big\{(\bar v^{ij} -U^{ij}(\vec{\bar  v}))(t,x); \min\Big\{(\bar  v^{ij }-L^{ij}(\vec{\bar v}))(t,x);\\\qq   -\partial_t \bar v^{ij}(t,x) -\mathcal{L}^{X}(\bar v^{ij})(t,x) -f^{ij}(t,x, (\bar v^{kl}(t,x))_{\kl},\sigma(t,x)^\top D_x \bar v^{ij}(t,x)) \Big\}  \Big\}  =0,\\
  \bar  v^{ij}(T,x)=h^{ij}(x)). \\
           \end{array}          \right. \end{eqnarray}
Thus $(\bar v^{ij})_{\ij}$ is also a solution for the multi-dimensional system (\ref{system2}) and then,
by uniqueness of the solution of (\ref{system2}) in $\Pi_g$, we have $\bar v^{ij}=\udl v^{ij}$ for any $\ij$, which completes the proof. 

\begin{rem} \label{L2-estimates}     The result of Theorem \ref{thmprincipalsection1} is still valid if (H0)-(H4) are in
force and the functions $(\underline g_{ij})_{\ij}$ verify (H5) since, by symmetry, one can go through the decreasing scheme (\ref{penalizedscheme}) to the increasing one (\ref{second-penalizedscheme}) and conversely. \qed
\end{rem}

Next, let us introduce the following family of processes $(Y^{ij})_{\ij}$ defined through the common solution $(v^{ij})_{\ij}$ of the min-max and max-min systems as follows: $\forall s\leq T$ and $\ij$, 
\be \label{Feynmankac-rep} 
Y^{ij}_s=v^{ij}(s\vee t, X^{t,x}_{s\vee t}).
\ee
We are going to show that the backward SDE counterpart of system (\ref{system1}) (or (\ref{system2})) has a unique global solution. Actually, we have:   
\begin{thm}\label{existence-globalsol}Assume that Assumptions (H0)-(H5) are fulfilled. Then there exist processes 
$(Z^{ij})_{\ij}$, $(K^{ij,+})_{\ij}$ and 
$(K^{ij,-})_{\ij}$ which belong respectively to $\cH^{2,d}_t$, $ \cA^{2}_t $ and $ \cA^{2}_t $ (which depend on $(t,x)$ and which we omit to precise) such that  
the family $(Y^{ij}, \; Z^{ij},\; K^{ij,+},\; K^{ij,-})_{\ij}$ is a solution of the following doubly reflected BSDEs (DRBSDE in short) with bilateral interconnected obstacles: For any $\ij$ and $s\in [t,T]$, 
 \be \label{dreflect-bsde} \left\{ 
 \begin{array}{l}  
\displaystyle{ d Y_{s}^{ij} = -f^{ij}(s,X^{t,x}_s, \vec{Y}_s, Z_{s}^{ij}) ds + dK_{s}^{ij,-} -dK_{s}^{ij,+} -Z_s^{ij}dB_s \,\,;\, Y^{ij}_T=h^{ij}(X^{t,x}_T);}
 \\\\
\displaystyle{ Y_s^{ij} \le U^{ij}_s(\vec{Y})  \; \textrm{and} \; Y_s^{ij} \ge
  L^{ij}_s(\vec{Y})\,\,;
} \\\\
\int_{t}^{T} (Y_s^{ij} - U^{ij}_s(\vec{Y})) dK_{s}^{ij, -} =0\;\textrm{and} \;\int_{t}^{T} (  L^{ij}_s(\vec{Y}) - Y_s^{ij}
) dK_{s}^{ij, +} =0\\ 
 \end{array}  \right. \ee
where, for each $\ij$, the lower (resp. upper)
 interconnected obstacle $ L^{ij}(\vec{Y})$ (resp. $ U^{ij}(\vec{Y})$)
is defined by: $\forall s\in [t,T]$, $$\displaystyle{ L_s^{ij}(\vec{Y}) = \max_{k \in \gami} \{ Y^{kj}_{s} - \underline{g}_{ik}(s,X_s^{t,x})\} \; (\textrm{resp.} \; 
 U_{s}^{ij}(\vec{Y}) =  \min_{l \in \gamj} \{ Y^{il}_{s} + \bar{g}_{jl}(s,X_s^{t,x})\}). }$$

This solution is unique in the following sense: If $(\bar Y^{ij}, \bar Z^{ij},\bar K^{ij,+},\bar K^{ij,-})_{\ij}$ is another solution of (\ref{dreflect-bsde}) then for any $(i,j)\in \Gamma$, $Y^{ij}=\bar Y^{ij}$,  $Z^{ij}=\bar Z^{ij}$ and 
$K^{ij,+}-K^{ij,-}=\bar K^{ij,+}-\bar K^{ij,-}$.  
 
\end{thm}
\no \udl{{\bf Proof}}: It is postponed to Appendix (A1) relegated to the end of this paper.
\section{The min-max (or max-min) solution as the value of the zero-sum switching game}
In this section, our objective is to study the link of the solution $(v^{ij})_{\ij} $ of both the min-max and max-min system with the values of an explicit switching game. 
To do this, we shall deeply rely on the representation (\ref{Feynmankac-rep}) in terms of the solution $(Y^{ij})_{\ij}$ of the general DRBSDE given in Theorem \ref{existence-globalsol}. 

So once for all in this section, we suppose that Assumptions (H0)-(H5) hold. On the other hand we assume that: 

\noindent {\bf (H6)}:
$$\begin{array}{l}\mbox{ For any }\ij, \mbox{ the function } f^{ij}\mbox{  does not depend on }(\vec y,z^{ij}). \qq\qq\qq
\end{array}
$$
\subsection{Description of the zero-sum switching game}Assume we have two players $\pi_1$ and $\pi_2$ who intervene on a
system (e.g. the production of energy from several sources such as oil, cole, hydro-electric, etc.) with the help of switching strategies. An admissible switching strategy
for $\pi_1$ (resp. $\pi_2$) is a sequence $\d:=(\sigma_n,\xi_n)_{n\ge 0}$
(resp. $\nu:=(\tau_n,\zeta_n)_{n\ge 0}$) where 
for any $n\geq 0$,

(i) $\sigma_n$ (resp. $\tau_n$) is an $\bF$-stopping times such that
$\P$-a.s., $\sigma_n \leq \sigma_{n+1}\leq T$ (resp. $\tau_n\leq \tau_{n+1}\leq T$) ;

(ii) $\xi_n$ (resp. $\zeta_n$) is a random variable with values in $\Gamma^1$
(resp. $\Gamma^2$) which is $\cF_{\sigma_n}$ (resp.
$\cF_{\tau_n}$)-measurable ;

(iii) $\P[\sigma_n<T, \forall n\geq 0]=\P[\tau_n<T, \forall n\geq
0]=0$ ;

(iv) If $(A^\d_s)_{s\le T}$ and $(B^\nu_s)_{s\leq T}$ are the $\bF$-adapted RCLL processes defined by:
$$\forall \; s\in [t,T),\quad A^\d_s=\sum_{n\geq
1}\underline {g}_{\xi_{n-1}\xi_n}(\sigma_n,
X^{t,x}_{\sigma_n})1_{[\sigma_n\leq s]}\quad \mbox{ and
}\quad A^\d_T=\lim_{s\rw T}A^\d_s,
$$
and
$$\forall \; s\in [t,T), \quad B^\n_s=\sum_{n\geq
1}\bar{g}_{\z_{n-1}\z_n}(\t_n, X^{t,x}_{\t_n})1_{[\t_n\leq s]}\quad
\mbox{ and }\quad B^\n_T=\lim_{s\rw T}B^\n_s
$$
then $\E[(A^\d_T)^2+(B^\nu_T)^2]<\infty$.
For any $s\le T$, $A^\d_s$ (resp. $B^\nu_s$) is the cumulative switching cost at time $s$ for $\pi_1$ (resp. $\pi_2$) when she implements the strategy $\d$ (resp. $\nu$).\ms

Next let $\ig$  and $t\in [0,T]$ be fixed. We say that the admissible switching strategy $\d:=(\sigma_n,\xi_n)_{n\ge 0}$
(resp. $\nu:=(\tau_n,\zeta_n)_{n\ge 0}$) of $\pi_1$ (resp. $\pi_2$) belongs to $\cA^i_{\pi_1}(t)$ (resp. $\cA^j_{\pi_2}(t)$) if $\sigma_0=t,\,\, \xi_0=i$ (resp. $\t_0=t,\,\, \z_0=j$).

Given an admissible strategy $\d$ (resp. $\n$) of $\pi_1$ (resp. $\pi_2$), one associates
a stochastic process $(a_s)_{s \le T}$ (resp. $(b_s)_{s\leq T}$) which indicates along with time the current mode of $\pi_1$ (resp. $\pi_2$) and which is defined by: 
\be \label{defuv}\forall s \le T,\,\,
a_s=\xi_0 1_{\{\sigma_0\}}(s)+\sum_{n\geq 1}\xi_{n-1}
1_{]\sigma_{n-1},\sigma_{n}]}(s) \,\,(\mbox{resp.}\,\,\, b_s=\zeta_0 1_{\{\t_0\}}(s)+\sum_{n\geq 1}\zeta_{n-1}
1_{]\t_{n-1},\t_{n}]}(s)).
\ee

Let now $\delta=(\sigma_n,\xi_n)_{n\ge 0}$ (resp.
$\nu=(\tau_n,\zeta_n)_{n\ge 0}$) be an admissible strategy for $\pi_1$ (resp. $\pi_2$) which belongs to ${\cal A}^i_{\pi_1}(t)$ (resp. ${\cal A}^j_{\pi_2}(t)$). The interventions of the players are not free and generate a payoff which is a reward (resp. cost) for
$\pi_1$ (resp. $\pi_2$) and whose expression is given by
\be \label{cost_functional}\begin{array}{l}
J_{t}(\d,\n):=\E[h^{a_{T} b_T}(X_T^{t,x})+\int_{t}^T f^{a_r b_r}(r,X^{t,x}_r)dr-A^\d_T+B^\nu_T|\cF_t].\ea
\ee
When the system is in $(i,j)$ at the initial time $\ot$, we define the upper (resp. lower) value of the game by $$ \bar V_t^{ij} :=
\essinf_{\nu\in {\cal A}^j_{\pi_2}(t)}\esssup_{\d\in {\cal A}^i_{\pi_1}(t)} J_{t}(\d,\nu) \,\,(\mbox{resp}.\, \underbar V_t^{ij}=\esssup_{\d\in {\cal A}^i_{\pi_1}(t)}\essinf_{\nu\in {\cal A}^j_{\pi_2}(t)}
J_{t}(\d,\nu) ).
$$ The game has a value if $\bar V_t^{ij}=\underbar V_t^{ij}$ and finally, we say that the game has a saddle-point if there exists a pair of admissible strategies $(\d^*,\nu^*)\in 
{\cal A}^i_{\pi_1}(t)\times {\cal A}^j_{\pi_2}(t)$ such that for any $\d$ and $\nu$, it holds
\begin{equation}\label{saddlepoint}
J_{t}(\d,\nu^*)\le J_{t}(\d^*,\nu^*)\le J_{t}(\d^*,\nu).
\end{equation}

As previously mentionned, we are going to study the link between the solution $\vij$ of (\ref{system1}) and the upper and lower values   
$\bar V_t^{ij}$, $\underbar V_t^{ij}$ of the zero-sum switching game.
\begin{thm}\label{gamevalue}
Suppose that Assumptions (H0)-(H6) are in force. Then the processes $(Y^{ij})_{\ij}$ of the unique solution of the doubly reflected BSDE (\ref{dreflect-bsde}) satisfy
\be\label{lienyijvjeux}\underbar{V}_t^{ij}  \le Y_{t}^{ij}=v^{ij}(t,x) \le \bar V_t^{ij}. 
\ee
\end{thm}
\no \udl{{\bf Proof}:} Let us consider two families of auxiliary processes $ (\hat{U}^{\delta j})_{j\in \Gamma^2}$ and
$(\hat{U}^{i \nu})_{i\in \Gamma^1}$ associated 
 with admissible strategy $\delta \in \cA^i_{\pi_1}(t)$ and $\nu \in \cA^j_{\pi_2}(t)$ and defined by: $\forall j\in \gmj$, 
 \be\label{equationhatudelta}
  \left\{ \begin{array}{ll}
  &\hat U^{\d j}\in \cS^2_{t, d},\,\hat Z^{\d j}\in \cH^{ 2,d} _t,\,\, \hat K^{\delta
 j, -}\in \cA^2_{t, i};\\\\
 &\hat{U}_{s}^{\delta j}=\; \; 
 h^{a_T j}(X^{t,x}_T) +\int_{s}^{T} f^{ a_r j}(r,X^{t,x}_r) dr -\int_{s}^{T} \hat Z_{r}^{\delta j}dB_{r} -(A_{T}^{\delta} -A_{r}^{\delta}) -
 (\hat K^{\delta j, -}_{T} -\hat K_{s}^{ \delta
 j, -}),\,\,s\in [t, T]; \\\\
&\forall s\in [t, T], \, \hat{U}_{s}^{\delta j} \le \min_{l \neq j}
 \left(\hat{U}_{s}^{\delta  l} +\bar{g}_{jl}(s,X^{t,x}_s)\right) \mbox{ and }
 \int_{t}^{T}\{\hat{U}_{r}^{\delta  j} -\min_{l \neq j}
\{\hat{U}_{r}^{\delta  l} +\bar{g}_{jl}(r,X^{t,x}_r)\}\}d\hat K_r^{\delta j, -} =0.
\end{array} \right.
 \ee
and for any $i\in \gmi$
\be \label{equationunu}
  \left\{ \begin{array}{ll}
 {}&U^{i \nu}\in \cS^2_{ t,d},\, Z^{i \nu}\in \cH^{ 2,d}_t,\,\, K^{i \nu, +}\in \cA^2_{ t,i};\\  \\
 {}& U_{s}^{i \nu}=\;  h^{i b_T}(X^{t,x}_T)  +\int_{s}^{T}f^{ i b_r}(r,X^{t,x}_r)dr -\int_{s}^{T}Z_{r}^{i \nu}dB_{r} 
 +(B_{T}^{\nu} -B_{s}^{\nu}) +(K^{i \nu, +
 }_{T} -K_{s}^{i \nu, , +}),\,s \in [t,T] \,; \\\\ {}&\forall s\in [t,T],\,\,
 U_{s}^{i \nu} \ge \max_{k \neq i}
\{U_{s}^{k \nu} - \underline{g}_{ik}(s,X^{t,x}_s)\}\mbox{ and }\int_{t}^{T}\left( \hat{U}_{r}^{\delta j} -\max_{k \neq i}
\{U_{r}^{k \nu} - \underline{g}_{ik}(r,X^{t,x}_r)\}\right) dK_r^{ i \nu,+ } =0.
 \end{array} \right.
 \ee

\no These equations are actually not of standard form, but once more by a change of variables, one verifies that both $ (\hat{U}^{\delta j} - A^{\delta})_{j \in \Gamma^2}$ and $(U^{i \nu} +B^{\nu})_{i \in \Gamma^1} $ solve standard multi-dimensional RBSDEs system which have a unique solution. On the other hand, let us point out that thanks to the connection between
the standard switching problem and multi-dimensional RBSDE with upper (resp. lower) interconnected obstacles (see e.g. \cite{Djehichehampopier} or \cite{Hutang}) the family 
$(\hat{U}^{\delta j}-A^{\delta})_{j \in \Gamma^{2}}$ (resp. $(U^{i \nu }+B^{\nu})_{i \in \Gamma^1}$) of processes verifies:
$$ \hat{U}_{t}^{\delta j}-A^\d_t = \essinf_{\nu \in \mathcal{A}^j_{\pi_2}(t)}J_{t}(\delta, \nu)\,\,
\textrm{(resp.  } \;  U_{t}^{i \nu }+B^\nu_t =\esssup_{\delta  \in \mathcal{A}^i_{\pi_1}(t)} J_{t}(\delta, \nu))$$
and then 
\be\label{equaiuessinf} \hat{U}_{t}^{\delta j} = \essinf_{\nu \in \mathcal{A}^j_{\pi_2}(t)}\{J_{t}(\delta, \nu)+A^\d_t\}\,\,
\textrm{(resp.  } \;  U_{t}^{i \nu } =\esssup_{\delta  \in \mathcal{A}^i_{\pi_1}(t)} \{J_{t}(\delta, \nu)-B^\nu_t\}).\ee 
  In order to prove (\ref{lienyijvjeux}), it is enough to establish the following:
   \begin{equation}\label{lowerupperctrl-drbsdesol}   
\forall \; \d \in \mathcal{A}^i_{\pi_1}(t),\; \nu \in \mathcal{A}^j_{\pi_2}(t)\;\;  \hat{U}_{t}^{\delta j} -A_t^{\d} \le Y_t^{ij} \le    
 U_{t}^{i \nu } +B_t^{\nu}, \end{equation}
   which equivalently means that 
	$$ \esssup_{\{\delta  \in \mathcal{A}^i_{\pi_1}(t)\} } \{ \hat{U}_{t}^{\delta j} -A_t^{\d} \} \le Y_t^{ij} \le    
\essinf_{ \{\nu \in \mathcal{A}^j_{\pi_2}(t) \}} \{ U_{t}^{i \nu } +B_t^{\nu} \},$$and the result follows by (\ref{equaiuessinf}). 

In the sequel, we only prove the first inequality in (\ref{lowerupperctrl-drbsdesol}) since the second one can be obtained by symmetry comparing $Y_t^{ij}$ to $U_t^{i \nu} + B_t^{\nu}$ for an arbitrary $\nu$.

So let us consider, on the time interval $[t,T]$, the penalized decreasing
scheme introduced in (\ref{intermpenscheme}). The processes \\ $( Y^{ij,m},  Z^{ij,m}, K^{ij,m})_{\ij}$, $m\geq 0$, verify: $\forall \ij$, \be\label{intermpenschemebis}\left\{\begin{array}{l} 
Y^{ij,m}\in \mathcal{S}^2_{t}, \,\, Z^{ij,m} \in \mathcal{H}_t^{2,d} \mbox{ and
} K^{ij,m } \in \mathcal{A}^{2}_{t, i}\,\,;\\
 Y^{ij,m}_s=h^{ij}(X^{t,x}_T)+\int_{s}^{T}
f^{ij,m}(r,X^{t,x}_r,(Y^{kl,m}_r)_{\kl})dr + ( K^{ij,m}_T-
 K^{ij,m}_s)-\int_s^T Z^{ij,m}_rdB_r,\,\,\forall s\in [t,T];\\
\displaystyle{ Y^{ij,m}_{s}\geq \max_{k\in
(\Gamma^{2})^{-i}}\{ Y^{kj,m}_s -\underline{g}_{ik}(s,X^{t,x}_s)\}},
\,\,\forall s\in [t,T];\\
\int_t^T( Y^{ij,m}_s - \max_{k\in
(\Gamma^{1})^{-i}}\{ Y^{kj,m}_s -\underline{g}_{ik}(s,X^{t,x}_s)\})d
 K^{ij,m}_s=0\end{array}\right.\ee
where, we recall that
$$f^{ij,m}(s, X^{t,x}_s,\vec{y}) = f^{ij}(s,X^{t,x}_s) - m \sum_{l \in \gamj} \big(y^{ij} -(y^{il} + \bar
{g}_{jl}(s,X^{t,x}_s)) \big)^{+}.
$$
As already mentioned, we know that, for any $\ij$, $Y^{ij,m}\rw_m Y^{ij}$ in $\cS_t^2$. \\

Next fix $(\ijo)\in \gam$ and let us show that $Y^{i_0j_0}_t \ge \hat{U}_t^{\delta j_0}-A^\d_t$ for any $\d=(\sigma_l,\xi_l)_{l\ge 0}$ in ${\cal A}^{\io}_{\pi_1}(t)$.
 So let us define the processes $(Y^{\delta j, m})_{j\in \gmj}$  and $(\hat{U}^{\delta j, m})_{j\in \gmj}$ as follows: 
\ms

\noindent (i) $\forall j\in \Gamma^2$, $$\forall s\in [t,T),\quad Y^{\delta j, m}_s=\sum_{l \geq 0}Y^{\xi_l j, m}_s\ind_{[\sigma_l\le s<\sigma_{l+1}]}\quad \mb{ and
}\quad Y^{\delta j, m}_T=h^{a_T j}(X^{t,x}_T),
$$ where,
\be\label{ystrategy}
\forall s\in [t,T],\quad Y^{\xi_lj,m}_s=\sum_{q\in
\Gamma^1} Y^{qj, m}_s\ind_{[\xi_l=q]}.
\ee
The process $Y^{\d j,m}$ is well defined since the sum contains
only  finitely many terms as the strategy $\d$ is admissible and then $\P[\si_l<T, \forall l\geq 0]=0$. On the other hand, at time $0<\sigma_l<T$, $Y^{\delta j, m}$ has a jump which is
equal to $Y^{\xi_lj, m}_{\sigma_l}-Y^{\xi_{l-1}j,m}_{\sigma_l}.$ 
\ms

\noindent (ii) The processes $(\hat{U}^{\delta j ,m})_{j\in \gmj}$ 
are defined as the solution in  $\cS_{t, d}^2$ of the following non standard multi-dimensional BSDE: $\forall j\in \gmj$,  
\be\label{penalizedhatudelta}
  \begin{array}{ll}
\hat{U}_{s}^{\delta j,m}=\; &\; h^{a_T j}(X^{t,x}_T) +\int_{s}^{T} \Big \{f(r,X^{t,x}_r, a_r, j)
-m \sum_{l \ne j}(\hat{U}_r ^{\delta j, m} -\hat{U}_r ^{\delta l,m} -\bar{g}_{jl} )^{+}\Big \}dr \\\\
 \;&\;\qq\qq -(A_{T}^{\delta} -A_{s}^{\delta})-\int_{s}^{T}\hat V_{u}^{\delta j,m}dB_{u},\,\,s\in [t,T]. 
\end{array}
\ee
Note that $(\hat{U}^{\delta j ,m}+A^\d)_{j\in \gmj}$ is a solution of a standard multi-dimensional BSDE whose coefficient is Lipschitz. As those latter processes exist, then so are $(\hat{U}^{\delta,j ,m})_{j\in \gmj}$. 
On the other hand, as for the system given in (\ref{intermpenschemebis}), 
the sequence of processes $((\hat{U}^{\delta j ,m})_{j\in \gmj})_{m\ge 0}$ converges in $\cS_{t, d}^2$ toward
$(\hat{U}^{\delta j})_{j\in \gmj}$.\\ We now prove the following: for any $m\ge 0$ and $j\in \gmj$,
\begin{equation}
 \label{auxiliar_ineq}  Y^{\delta j,m}_t\ge \hat{U}_t^{\delta j,m}.
\end{equation}
For any $j\in \gmj$,
let us define $K^{\delta j, m}$ and $Z^{\delta j, m}$ as follows:
$\forall s\in [t,T]$,
$$Z^{\delta j, m}_s:=\sum_{l\geq 0}
Z^{\xi_l j, m}_s\ind_{[\sigma_l\le s <\sigma_{l+1}[}\quad 
\mbox{ and }\quad K^{\delta j, m}_s=\sum_{l\geq 0}\int_{s\wedge \sigma_l}^{s\wedge\sigma_{l+1}}
dK^{\xi_lj, m}_s,
$$
where, $Z^{\xi_lj, m}_s$ and $K_s^{\xi_lj, m}$ are defined in the same way as in (\ref{ystrategy}). Once more, there is no definition issue of those processes since $\d$ is admissible. Therefore the triple of processes $(Y^{\delta j, m},Z^{\delta j, m},K^{\delta j, m})_{j\in \gmj}$ 
verifies: $\forall s\in [t,T)$, 
$$\begin{array}{ll}
Y^{\delta j, m}_s&=Y^{\delta j, m}_t-\int_t^{s}\Big \{f^{a_r j}(r,X^{t,x}_r)dr
+m\sum_{l\neq j}\big(Y^{\delta j, m}_r-Y^{\delta l, m}_r - 
\bar{g}_{jl}(r,\x_r)\big)^{+}dr+
Z^{\delta j, m}_rdB_r-dK^{\delta j, m}_r\Big \}\\&\qq\qq\qq+
\sum_{l\geq 1}
(Y^{\xi_lj, m}_{\sigma_l}-Y^{\xi_{l-1}j, m}_{\sigma_l})\ind_{[\sigma_l\leq s]}\\
{}&=Y^{\delta j, m}_t-\int_t^{s}\Big \{f^{a_rj}(r,X^{t,x}_r)dr
+m\sum_{l\neq j}\big(Y^{\delta,j, m}_r-Y^{\delta l, m}_r - 
\bar{g}_{jl}(r,\x_r)\big)^{+}dr+
Z^{\delta j, m}_rdB_r-dK^{\delta j, m}_r\Big \}\\&\qq\qq\qq-
\sum_{l\geq 1}
(Y^{\xi_{l-1}j, m}_{\sigma_l}-Y^{\xi_lj, m}_{\sigma_l}+\underline{g}_{\xi_{l-1}\xi_{l}}(\sigma_l, X^{t,x}_{\sigma_l}))\ind_{[\sigma_l\leq s]}+A^\d_s. 
\end{array}
$$
Next, let us define $\tilde A^{\d j,m}$ by: 
\be \label{def-tildAproc}
\tilde A^{\d j,m}_s:=\sum_{l\geq 1}
(Y^{\xi_{l-1}j, m}_{\sigma_{l}}-Y^{\xi_{l}j, m}_{\sigma_l}+\underline{g}_{\xi_{l-1} \xi_{l}}(\sigma_l, X^{t,x}_{\sigma_l}))
\ind_{[\sigma_l\leq s]}\;\;\mbox{ for }\;\; s\in [t,T)\;\; \mbox{ and }\;\;  \tilde A^{\d j}_T=\lim_{s\rw T}
\tilde A^{\d j}_s,\ee which is an adapted non-decreasing process. As the strategy $\d$ is admissible, then writing backwardly
between $s$ and $T$ the equation for the process $Y^{\d j, m}$ 
we obtain: $\forall j\in \Gamma^2$, 
\be\label{eqydjm}\begin{array}{ll}
Y^{\delta j, m}_s&=h^{a_T j}(X^{t,x}_T)+\int_s^{T}\Big \{f^{a_r j}(r,X^{t,x}_r)dr
-m\sum_{l\neq j}\big(Y^{\delta j, m}_r-Y^{\delta l, m}_r - 
\bar{g}_{jl}(r,\x_r)\big)^{+}dr\\{}&\qq\qq-
Z^{\delta j, m}_rdB_r+dK^{\delta j, m}_r\Big \}-(A^\d_T-A^\d_s)+(\tilde A^{\d j,m}_T-\tilde A^{\d j,m}_s),\,\forall s \in [t,T]. 
\end{array}
\ee
This equation implies also that $\E[(\tilde A^{\d j,m}_T)^2]<\infty$. 

Let us now introduce the sequence of processes  
$((Y^{\delta j, m,k},Z^{\delta j, m,k})_{j\in \G^2})_{k\geq 0}$ ($m$ is fixed) defined by Picard type iterations as follows:

\no (a) For any $j\in \G^2$ and $s\in [t,T]$, $(Y^{\delta j, m,0}_s,Z^{\delta j, m,0}_s)=(\hat U^{\d j,m}_s,
\hat V^{\d j,m}_s)$ ;

\no (b) For $k\geq 1$, $(Y^{\delta j, m,k},Z^{\delta j, m,k})_{j\in \G^2}$ is the solution in $\cS^2_{t,d}\times \cH^{2,d}_t$ of the following multidimensional BSDE: $\forall j\in \G^2$
\be\label{eqydjmk}\begin{array}{ll}
Y^{\delta j, m,k}_s&=h^{a_T j}(X^{t,x}_T)+\int_s^{T}\Big \{f^{a_r j}(r,X^{t,x}_r)dr
-m\sum_{l\neq j}\big(Y^{\delta j, m,k}_r-Y^{\delta l, m,k-1}_r - 
\bar{g}_{jl}(r,\x_r)\big)^{+}dr\\{}&\qq\qq-
Z^{\delta j, m,k}_rdB_r\Big \}-(A^\d_T-A^\d_s)+(K^{\delta j, m}_T-K^{\delta j, m}_s)+(\tilde A^{\d j,m}_T-\tilde A^{\d j,m}_s),\,\forall s \in [t,T] 
\end{array}
\ee
First note that those backward equations are decoupled and the processes $(Y^{\delta j, m,k},Z^{\delta j, m,k})$ are well defined by an obvious change of variables.  Next by the one dimensional comparison theorem \cite{epq} one has $Y^{\delta j, m,1}\geq \hat U^{\d j,m}$ and $Y^{\delta j, m,k}\geq Y^{\delta j, m,k-1}$, $k\ge 1$ and $j\in \G^2$. This last inequality is obtained by induction since for any $q \neq j$, the mapping 
$y_q \in \R \mapsto 
-m\sum_{l\neq j}\big(y_j-y_l - 
\bar{g}_{jl}(r,\x_r)\big)^{+}$ is non-decreasing when the other variable are fixed. Finally since the coefficients of the BSDEs of (\ref{eqydjmk}) are Lipschitz then for any $j\in \G^2$, the sequence 
$(Y^{\delta j, m,k},Z^{\delta j, m,k})_{k\ge 0}$ converges in $\cS^2_{t,d}\times \cH^{2,d}_t$ to a limit 
$(\bar Y^{\delta j, m},\bar Z^{\delta j, m})$ which satisfies: $\forall j\in \G^2$, 
\be\label{eqydjmktilde}\begin{array}{ll}
\bar Y^{\delta j, m}_s&=h^{a_T j}(X^{t,x}_T)+\int_s^{T}\Big \{f^{a_r j}(r,X^{t,x}_r)dr
-m\sum_{l\neq j}\big(\bar  Y^{\delta j, m}_r-\bar  Y^{\delta l, m}_r - 
\bar{g}_{jl}(r,\x_r)\big)^{+}dr\\{}&\qq\qq-
\bar Z^{\delta j, m}_rdB_r\Big \}-(A^\d_T-A^\d_s)+(K^{\delta j, m}_T-K^{\delta j, m}_s)+(\tilde A^{\d j,m}_T-\tilde A^{\d j,m}_s),\,\forall s \in [t,T] 
\end{array}
\ee
But the solution of (\ref{eqydjm}) is unique then for any $j\in \G^2$, 
$$
Y^{\delta j, m}=\bar Y^{\delta j, m}=\lim_{k\rw \infty}Y^{\delta j, m,k}\geq \hat U^{\d j,m}
$$
which ends the proof of the claim. 

Taking now the limit w.r.t. $m$, we obtain that 
$$
Y_{t}^{i j}=\lim_{m\rw \infty} Y^{i j,m}_t=\lim_{m\rw \infty} \{Y_{t}^{\delta j, m}-A^\d_t\} 
\geq \lim_{m\rw \infty} \{\hat U^{\d j,m}_t-A^\d_t\}=\hat U^{\d j}_t-A^\d_t,\quad 
\forall j\in \gmj.
$$
which therefore yields
\begin{equation}\label{lowervalueineq} 
Y_{t}^{i j} \ge \esssup_{\d\in \cA^i_{\pi_1}(t)} (\hat U^{\d j}_t-A^\d_t) =
\underbrace{\esssup_{\d} \essinf_{\nu} J_t(\d, \nu)}_{ := \underbar{V}_t}. 
\end{equation}
Similarly as above and introducing first the sequence ($U^{i\nu, n}$)$_{n \ge 0}$ such that, for any $n$, $U^{i \nu,n}$
is obtained from $ U^{i \nu} $ (defined in (\ref{equationunu})) 
by penalization of the lower barrier, one can prove that 
$$ \udl Y_{s}^{i \nu , n}\leq  U^{i \nu,n}_s\quad \forall s\in [t,T],$$
 by comparing the two penalized equations
solved by $ \udl{Y}^{ij, n}$ and $U^{i \nu, n}$.
 Next, relying on the fact that $ ( \udl{Y}^{ij, n})$ is related to the increasing penalized scheme (converging to $Y^{ij}$)
one gets
$$Y_t^{ij} =\lim_{n} \udl{Y}^{ij, n}_t = \lim_{n}\{ \udl{Y}^{i\n, n}_t   +B_t^{\nu}\} \le
\lim_{n}\{ U_t^{i \nu, n} +B_t^{\nu}\} := U_t^{i \nu} +B_t^{\nu}. $$
Taking the infimum over all admissible strategies $\nu$ and reminding the interpretation of $
  U^{i \nu} +B^{\nu},$ one obtains 
 \begin{equation}\label{uppervalue-ineq}
Y_{t}^{i j} \le \essinf_{\nu \in \cA^j_{\pi_2}(t)} ( U^{i \nu }_t+ B^{\nu}_t) =
\underbrace{ \essinf_{\nu \in \cA^j_{\pi_2}(t)} \esssup_{\d \in \cA^i_{\pi_1}(t)}  J_t(\d, \nu)}_{ := \overline{V}_t}. 
\end{equation}  
 which achieves the proof.  \qed
\ms

The point now is which kind of additional assumptions should we add in order to have equalities in (\ref{lienyijvjeux}) and then the game has a value. The response to this question is affirmative if we moreover assume that the utilities $f_{ij}$ and $h_{ij}$ are separated with respect to $i$ and $j$, i.e., if they satisfy the following assumption: 
\ms

\noindent {\bf (H7)}: the two families ($f^{ij}$)$_{\ij}$ and ($ h^{ij}$)$_{\ij}$ of functions satisfy, for any $\ij$, 
$$   \quad f^{ij} : = f_{1}^{i} + f_{2}^{j} \;\mathrm{and}\;  h^{ij} : = h_{1}^{i} + h_{2}^{j}.$$

Once more we recall that we work under the assumptions (H0)-(H7). So let us consider the following system of reflected BSDEs with one interconnected lower (resp. upper) obstacles associated with \\$((f_{1}^{i})_{i\in \Gamma_1}, (h_{1}^{i})_{i\in \Gamma_1}, (\underline g_{ik})_{i,k\in \Gamma_1})$ (resp. $((f_{2}^{j})_{i\in \Gamma_2}, (h_{2}^{j})_{j\in \Gamma_2}, (\bar g_{jl})_{j,l\in \Gamma_2})$): $\forall i\in \Gamma_1$ (resp. $j\in \Gamma_2$)
\be \left\{ 
\ba{l} \label{upperswitchinggame} 
\displaystyle{ d Y_{s}^{1,i} = -f_1^{i}(s,X^{t,x}_s) ds  -dK_{s}^{1,i} -Z_s^{1,i}dB_s, \,\,s\in [t,T]\,\mbox{ and }
 Y_{T}^{1,i}=h_1^{i}(X_T)};\\
\\
Y_s^{1,i} \ge
  L^{1,i}_s(\vec{Y}) := \max_{k \in \gami} ( Y^{1, k}_{s} - \underline{g}_{ik}(s,X_s^{t,x})),\,s\in [t,T]\,\textrm{and}
  \;\; \int_{t}^{T} (  L^{1,i}_s(\vec{Y}) - Y_s^{ij}
) dK_{s}^{1, i} =0 \\
\ea \right.
\ee
(resp. 
\be \left\{ 
\ba{l} \label{lowerswitchinggame}
 d Y_{s}^{2,j} = -f_2^{j}(s,X^{t,x}_s) ds  +dK_{s}^{2,j} -Z_s^{2,j}dB_s, s\in [t,T] \mbox{ and } Y_{T}^{2,j}=h_2^{j}(X_T);\\
\\
 \; Y_s^{2, j} \le
  U^{2,j}_s(\vec{Y}) := \min_{l \in \gamj} \left( Y^{2, k}_{s} + \bar{g}_{jl}(s,X_s^{t,x})\right),\,\,s\in [t,T]\,\,\textrm{and} \;\;
  \int_{t}^{T} ( Y_s^{2, j} -  U^{2,j}_s(\vec{Y})) dK_{s}^{2, j} =0).
 \\
\ea \right.
\ee 
Under assumptions (H3)-(H4), equations (\ref{upperswitchinggame}) and (\ref{lowerswitchinggame}) have unique solutions (see e.g. \cite{Hammorlais13}, Prop. 5.1, pp.188).

We then have:

\begin{prop}
\label{game-interp2} Suppose Assumptions (H0)-(H7) are 
satisfied. Then for any $\ij$, 
\ms

(i) $Y^{ij}=Y^{1,i}+Y^{2,j}$ ;

(ii) $\underbar{V}_t^{ij}=\bar V_t^{ij}=v^{ij}(t,x)$ and the game has a saddle-point $(\d^*,\nu^*)$.
\end{prop}
\no \pr: Let us first deal with the first point. For $\ij$, let $\underbar Y^{ij}=Y^{1,i}+Y^{2,j}$, $\underbar Z^{ij}=Z^{1,i}+Z^{2,j}$, $\underbar K^{ij,+}=K^{1,i}$ and $\underbar K^{ij,-}=K^{2,j}$. Then by (\ref{upperswitchinggame}) and (\ref{lowerswitchinggame}) it is easily seen that $(\underbar Y^{ij},\underbar Z^{ij},\underbar K^{ij,+},\underbar K^{ij,-})$ is also a solution for the DRBSE (\ref{dreflect-bsde}). As the solution of this latter is unique then $Y^{ij}=\underbar Y^{ij}=Y^{1,i}+Y^{2,j}$.

We now focus on the second point. First note that under the condition of separation of the utilities $f^{ij}$ and $h^{ij}$, for any 
$\d \in {\cal A}^i_{\pi_1}(t)$ and $\nu \in 
{\cal A}^j_{\pi_2}(t)$ we have 
\be \label{cost_functional}
J_{t}(\d,\n)=J^1_{t}(\d)+J^2_{t}(\n)
\ee where 
\be \label{cost_functional}\begin{array}{c}
J^1_{t}(\d):=\E[h_1^{a_{T} }(X_T^{t,x})+\int_{t}^T f_1^{a_r }(r,X^{t,x}_r)dr-A^\d_T|\cF_t] \mbox{ and }J^2_{t}(\nu):=\E[h_2^{b_{T} }(X_T^{t,x})+\int_{t}^T f_2^{b_r }(r,X^{t,x}_r)dr+B^\n_T|\cF_t].\end{array}
\ee Therefore we obviously have $$
\underbar{V}_t^{ij}=\bar V_t^{ij}=
sup_{\d \in {\cal A}^i_{\pi_1}(t)}J^1_t(\d)+\inf_{\n \in {\cal A}^j_{\pi_2}(t)}J^2_t(\n).
$$
and by (\ref{lienyijvjeux}) we obtain 
$\underbar{V}_t^{ij}=\bar V_t^{ij}=Y^{ij}_t=v^{ij}(t,x).$

Next the link between the solution of the BSDE (\ref{upperswitchinggame}) (resp. (\ref{lowerswitchinggame})) and the standard switching problem implies that (one can see \cite{Hammorlais13} or Step 2 of the proof of Theorem \ref{thmprincipalsection1} for more details):
$$Y^{1,i}_t=\sup_{\d \in {\cal A}^i_{\pi_1}(t)}J^1_t(\d)
\mbox{ and }Y^{2,j}_t=\inf_{\n \in {\cal A}^j_{\pi_2}(t)}J^2_t(\n).
$$
Finally let us define the strategies $\d^* = (\sigma_l^*,\xi_l^*)_{l \ge 0}$ \,\,(resp. $\nu^* =(\t_l^*,\z_l^*)_{l\geq 0} )$ as follows:
$$\sigma_0^*=t, \,\xi_0^*=i \mbox{ (resp. }\t_0^*=t,\; \z_0^*=j )$$ and for any  $l \ge 1$,
\be\label{optimstrat} \left\{
 \ba{l}  
 \sigma_{l}^* := \mbox{inf} \{ s \ge \sigma_{l-1}^*,\;  \;
  Y_s^{1, \xi^*_{l-1}} = 
 \mbox{max}_{k \neq \xi_{l-1}^*}
 \{Y_s^{1, k}  -\udl{g}_{ \xi_{l-1}^*, k}(s,X^{t,x}_s) \}  \} \wedge T \\
 \mbox{ and }\\
 \xi^*_{l}=\mbox{argmax}\{k \neq \xi_{l-1}^*, \,\, Y_{\sigma^*_l}^{1, k}  -\udl{g}_{ \xi_{l-1}^*, k}({\sigma^*_l},X^{t,x}_{\sigma^*_l}) \}
\ea
\right.\ee 
(resp.
 \be \label{optimstrat2} \left\{
 \ba{l}  
 \t_{l}^* := \mbox{inf} \{ s \ge \t_{l-1}^*,\;  \;
  Y_s^{2, \z^*_{l-1}} = 
 \mbox{min}_{q \neq \z_{l-1}^*}
 \{Y_s^{2, q}  +\bar{g}_{ \z_{l-1}^*, q}(s,X^{t,x}_s) \}  \} \wedge T \\
 \mbox{ and }\\
 \z^*_{l}=\mbox{argmin}\{q \neq \z_{l-1}^*, \,\, Y_{\t^*_l}^{2, q}  +\bar {g}_{ \z_{l-1}^*, q}({\t^*_l},X^{t,x}_{\t^*_l}) \}.
\ea
\right.\ee 
Then $\d^*$ (resp. $\n^*$) is admissible and optimal i.e., belongs to 
${\cal A}^i_{\pi_1}(t)$ (resp. ${\cal A}^j_{\pi_2}(t)$) and verifies 
 $$Y^{1,i}_t=\sup_{\d \in {\cal A}^i_{\pi_1}(t)}J^1_t(\d)=J^1_t(\d^*)
\mbox{ and }Y^{2,j}_t=\inf_{\n \in {\cal A}^j_{\pi_2}(t)}J^2_t(\n)=J^2_t(\n^*)
$$(see \cite{Hammorlais13} for more details).
Therefore 
$$J_{t}(\d,\n^*)=J^1_{t}(\d)+J^2_{t}(\n^*)\leq 
J^1_{t}(\d^*)+J^2_{t}(\n^*)=J_{t}(\d^*,\n^*)\leq 
J^1_{t}(\d^*)+J^2_{t}(\n^*)=J_{t}(\d^*,\n)
$$for any $\d \in {\cal A}^i_{\pi_1}(t)$ and $\n \in {\cal A}^j_{\pi_2}(t)$, which means that $(\d^*,\n^*)$ is a saddle-point for the zero-sum switching game. 
\begin{rem} Under (H0)-(H7), we also have the following relation: $\forall (i,j)\in \gam$, 
\begin{equation} \label{gamerepw1esp}\begin{array}{l}
v^{ij}(t,x)=
\sup_{\delta  \in \mathcal{A}^{i}_{\pi_1}(t)} \inf_{\nu  \in \mathcal{A}^{j}_{\pi_2}(t)}\E[\bar J_{t}(\delta, \nu) ]=\inf_{\nu  \in \mathcal{A}^{\jo}_{\pi_2}(t)}\sup_{\delta  \in \mathcal{A}^{\io}_{\pi_1}(t)} \E[\bar J_{t}(\delta, \nu) ].\end{array} \qed
\end{equation}
\end{rem}
\section{Conclusion}
In this paper, we have given appropriate conditions on the data
 of both the min-max and max-min systems so that their respective 
unique viscosity solutions coincide. These unique continuous viscosity solution have been constructed 
by means of a penalization procedure in the recent paper \cite{Djehichehammorlais14}. The main difficulty faced in that paper is that the two obstacles are interconnected and therefore not comparable. For this reason and without the comparison of the two barriers, we cannot apply the classical relationship between doubly reflected BSDEs, system of PDEs with lower and upper obstacles and the underlying game (see e.g. \cite{Hamhassani}). By providing appropriate regularity conditions so that comparison holds,  the viscosity solutions of the min-max and max-nin systems coincide. Finally, under further conditions on the drivers, this unique solution of the doubly reflected BSDE (\ref{dreflect-bsde}) is related and interpreted as the value function of a switching game.

\ms Since we also make use of the uniqueness of the viscosity solution of both the max-min and min-max systems to justify the uniqueness for the doubly reflected BSDE (\ref{dreflect-bsde}), our analysis deeply relies on the Markovian setting, therefore it seems quite natural to ask whether one can study the switching game 
in the general non-Markovian case. We leave this question for future research.
\section{Appendix}

\no {\bf (A1): Proof of Theorem \ref{existence-globalsol}}

Let us recall that for any $m\geq 0$, the quadruple $\yijm$ of (\ref{intermpenscheme}) verifies the following system of reflected BSDEs with one interconnected obstacle: $\forall s\in [0,T]$, 
\be\label{penscheme}
\left\{\begin{array}{l}
Y^{ij,m}\in \cS^{2}, \,\, Z^{ij,m} \in \cH^{ 2,d} \mbox{
and
} K^{ij,m,\pm } \in \cA_{i}^{2}\,\,;\\
 Y^{ij,m}_s=h^{ij}(X^{t,x}_T)+\int_{s}^{T}
{f}^{ij}(r,X^{t,x}_r,
( Y^{kl,m}_r)_{(k,l)\in \gam}, Z^{ij,m}_r)dr+\int_s^Td K^{ij,m,+}_r-\int_s^Td K^{ij,m,-}_r -\int_s^T Z^{ij,m}_rdB_r\,;\\
\displaystyle{ Y^{ij,m}_{s}\geq \max_{k\in
(\Gamma^{1})^{-i}}\{
Y^{kj,m}_s-\underline{g}_{ik}(s,X^{t,x}_s)\}}\,;\\
\int_0^T( Y^{ij,m}_s-\max_{k\in (\Gamma^{1})^{-i}}\{Y^{kj,m}_s-\underline{g}_{ik}(s,X^{t,x}_s)\})d
K^{ij,m}_s=0\end{array}\right.\ee 
where for any $s\in [0,T]$, $K^{ij,m,-}_s=m\int_0^s\{
\sum_{l\in \gamj}(Y^{ij,m}-Y^{il,m}_r-\bar g_{jl}(r,\xtx_r))^+\}dr$.

Let $(i,j)$ in $\Gamma_1 \times \Gamma^2$ be fixed. For $m\geq 0$ and $s\in [0, T]$, let us set
$ \alpha_s^{ij,m}:=\dfrac{dK_s^{ij, m-}}{ds}1_{[t,T]}(s)$. Then by estimate (\ref{eq:bound-penalizedterm2}) we have \be \label{estialpha}\begin{array}{l}
\E[\int_t^T| \alpha_s^{ij,m}|^2ds]\leq C_{tx},
\end{array}\ee
where $C_{tx}$ is a constant which may depend on $t$ and $x$. Therefore there exists a subsequence which we still denote by $\{m\}$ such that the sequence $(\alpha^{ij,m})_{m\ge 0}$ converges weakly to $\alpha^{ij}$ in $\cH^{2,1}.$ Also for $s\leq T$ let us set 
$$\ba{l}
k^{ij,-}_s=\int_t^{s\vee t}\a^{ij}_rdr.\ea
$$
Therefore the process $k^{ij,-}$ is continuous non-decreasing $\bF$-adapted and $\E[(k_T^{ij,-})^2]<\infty$. But by the 
representation property for any stopping time $\t \in [t,T]$, the sequence
$(\int_t^\t \alpha^{ij,m}_sds)_{m\ge 0}$ converges weakly in $L^2(\Omega, dP)$. 
\ms

Next, the sequence $(Y^{ij,m})_{m\geq 0}$ converges in $\cS^2_t$ to $Y^{ij}$. Actually, this stems from the uniform convergence of $(v^{ij,m})_{m\ge 0}$ to $v^{ij}$ in compact sets of $\espo$, the definition (\ref{Feynmankac-rep}) of $Y^{ij}$, estimate (\ref{estimationx}) and finally the polynomial growth of $v^{ij,m}$ and $\bar v^{ij}$ which comes from inequalities (\ref{inegvijnm}).  
\ms 

\no Now and using a classical method (see e.g. \cite{karoui}, 
proof of Theorem 5.2 in Section 6) and
since $(\a^{ij,m})_{m\ge 0}$ is uniformly bounded in $\cH^2_t$ then using It\^o's formula twice, respectively with $(Y^{ij,m})^2$ and $(Y^{ij,m}-Y^{ij,n})^2$, we deduce that:\\
(i) $$
\E[\int_t^T|Z^{ij,m}_s|^2ds+(K^{ij,m,+}_T)^2]\leq C_{t,x} \mbox{ and }$$
(ii) $$
\E[\int_t^T|Z^{ij,m}_s-Z^{ij,n}_s|^2ds]\rw_{n,m\rw \infty}0.$$
Thus, the sequence $(Z^{ij,m}1_{[t,T]}(.))_{m\ge 0}$ converges in 
$\cH_t^{2,d}$ to a limit which we denote by $Z^{ij}$. Finally and for any $s\in [t,T]$, let us set $$
k^{ij,+}_s=
-Y^{ij}_s+Y^{ij}_t-\int_t^s
{f}^{ij}(r,X^{t,x}_r,
( Y^{kl}_r)_{(k,l)\in \Gamma}, Z^{ij}_r)dr+ (k^{ij,-}_s-k^{ij,-}_t)-\int_t^sZ^{ij}_rdB_r\,\,(\mbox{note that } \;k^{ij,+}_t=0).$$
Therefore, for any stopping time $\t$, the subsequence 
$(K^{ij,m}_\tau-K^{ij,m}_t)_{m\ge 0}$ converges weakly in $L^2(dP)$ to $k^{ij,+}_\t$. Thus, the process $k^{ij,+}$ is non-decreasing and belongs to $\cS^2_t$. It follows that for any $s\in [t,T]$,\be\label{edsrlimite}
Y^{ij}_s= h^{ij}(X^{t,x}_T)+\int_s^T
{f}^{ij}(r,X^{t,x}_r,
( Y^{kl}_r)_{(k,l)\in \Gamma}, Z^{ij}_r)dr
+(k^{ij,+}_T-k^{ij,+}_s)-
(k^{ij,-}_T-k^{ij,-}_s )-\int_s^T Z^{ij}_rdB_r.\ee Henceforth the barriers $(\max_{k\neq i}(Y^{kj}_s-\underline g_{ik}(s,X^{t,x}_s)))_{s\in [t,T]}$ and 
$(\min_{l\neq j}(Y^{il}_s+\bar g_{jl}(s,X^{t,x}_s)))_{s\in [t,T]}$ verify the so-called Mokobodski assumption (see e.g. \cite{cvitanic-karatzas}), since $Y^{ij}$ is in between,
which is well-known in the two reflecting barriers BSDEs framework. Consequently, the double barrier reflected BSDE associated with
\\ $\{f^{ij}(s,X^{t,x}_s,(Y^{kl}_s)_{(k,l)\in \Gamma}, z),h^{ij}(X^{t,x}_T),
(\max_{k\neq i}(Y^{kj}_s-\underline g_{ik}(s,X^{t,x}_s)))_{s\in [t,T]}, (\min_{l\neq j}(Y^{il}_s+\bar g_{jl}(s,X^{t,x}_s)))_{s\in [t,T]}\}$
has a solution, i.e., there exist a quadruple of processes $(\tilde Y^{ij}_s,\tilde Z^{ij}_s, K^{ij,\pm}_s)_{s\in [t,T]}$ such that : $\forall s\in [t,T]$, 
 \be \label{dreflect-bsdebis} \left\{ 
 \begin{array}{l}  
\displaystyle{ d \tilde Y_{s}^{ij} = -f^{ij}(s,X^{t,x}_s, (Y^{kl}_s)_{(k,l)\in \Gamma}, \tilde Z_{s}^{ij}) ds + d\tilde K_{s}^{ij,-} -d\tilde K_{s}^{ij,+} -\tilde Z_s^{ij}dB_s \,\,; \tilde Y^{ij}_T=h^{ij}(X_T^{t,x});}
 \\
\tilde Y_s^{ij} \le U^{ij}_s(\vec{Y})\; \textrm{and} \; \tilde Y_s^{ij} \ge
  L^{ij}_s(\vec{Y})\,\,;
\\
 \int_{t}^{T} (\tilde  Y_s^{ij} - U^{ij}_s(\vec{Y})) d\tilde K_{s}^{ij, -} =0\;\textrm{and} \;\int_{t}^{T} (  L^{ij}_s(\vec{Y}) - \tilde Y_s^{ij}
) dK_{s}^{ij, +} =0\\ 
 \end{array}  \right. \ee (see \cite{cvitanic-karatzas} for more details).
On the other hand, since we are in the Markovian framework of randomness
(recall the representation (\ref{Feynmankac-rep}) for $(Y^{kl})_{(k,l)\in \Gamma}$), 
there exists a continuous deterministic function with polynomial growth such that 
$$
\forall s\in [t,T],\;\; \tilde Y^{ij}_s=\tilde v^{ij}(s,X^{t,x}_s). 
$$
Moreover the function $\tilde v^{ij}$  is the unique viscosity solution in $\Pi_g$ of the following PDE with two obstacles: 

\begin{eqnarray}\label{eqvijtilde} 
         \left\{ \begin{array}{l} 
   \min\Big\{(\tilde v^{ij }-L^{ij}(\vec {v}))(t,x),\max\Big\{ (\tilde  v^{ij} -U^{ij}(\vec { v}))(t,x), \\\qq\qq\qq-\partial_t \tilde  v^{ij} (t,x)-\mathcal{L}^{X}(\tilde v^{ij})(t,x) -f^{ij}(t,x, (v^{kl}(t,x))_{\kl},\sigma^\top(t,x)D_x\tilde  v^{ij}(t,x)) \Big\}  \Big\}  =0,\\
   \tilde  v^{ij}(T,x)=h^{ij}(x)\\
           \end{array}          \right. \end{eqnarray}(see \cite{Hamhassani}, pp. 261-262 or appendix (A2)). 
But $v^{ij}$ is also a solution of this latter, thus $v^{ij}=\tilde v^{ij}$ since the continuous solution
with polynomial growth of (\ref{eqvijtilde}) is unique. This uniqueness is due to polynomial growth of $f^{ij}(t,x,0,0)$ and $\vij$, and the continuity  of $(t,x)\mapsto f^{ij}(t,x,y,z)$ unifromly w.r.t. $(y,z)$. Then and for any $s\in [t,T]$, $Y^{ij}_s=\tilde Y^{ij}_s$.
Comparing now equations (\ref{edsrlimite}) and (\ref{dreflect-bsdebis}), we deduce that $\tilde Z^{ij}_s=Z^{ij}_s$ for any $s\in [t,T]$
which means that  the quadruple $(Y^{ij}_s, Z^{ij}_s, k^{ij,\pm}_s)_{s\in [t,T]}$ verifies: $\forall s\in [t,T]$, 
 \be \label{dreflect-bsdeeqxx} \left\{ 
 \begin{array}{l}  
\displaystyle{ d  Y_{s}^{ij} = -f^{ij}(s,X^{t,x}_s, (Y^{kl}_s)_{(k,l)\in \Gamma}, Z_{s}^{ij}) ds + dk_{s}^{ij,-} -dk_{s}^{ij,+} - Z_s^{ij}dB_s \,\,;\,Y^{ij}_T=h^{ij}(X_T^{t,x});}
 \\
\displaystyle{  Y_s^{ij} \le U^{ij}_s(\vec{Y})  \; \textrm{and} \; Y_s^{ij} \ge
  L^{ij}_s(\vec{Y})\,\,;
} \\
\int_{t}^{T} ( Y_s^{ij} - U^{ij}_s(\vec{Y})) d\tilde K_{s}^{ij, -} =0\;\textrm{and} \;\int_{t}^{T} (  L^{ij}_s(\vec{Y}) -  Y_s^{ij}
) dK_{s}^{ij, +} =0.\\ 
 \end{array}  \right. \ee But we can do the same for the other indices $(i_1,j_1) \in \Gamma$. Therefore, the processes 
 $((Y^{ij}_s, Z^{ij}_s, k^{ij,\pm}_s)_{s\in [t,T]})_{\ig}$ is a solution of system (\ref{dreflect-bsdebis}). \\

\noindent 
We now deal with the issue of uniqueness of the solution of (\ref{dreflect-bsdebis}). So let 
 $((\bar Y^{ij}_s, \bar Z^{ij}_s, \bar K^{ij,\pm}_s)_{s\in [t,T]})_{\ig}$ be another solution of (\ref{existence-globalsol}). We are going to show that for any $\ig$ and $m,n \geq 0$ we have, $\underbar Y^{ij,n}_s\le \bar Y^{ij}_s\leq \bar Y^{ij,m}_s$, $s\in [t,T]$, where $\bar Y^{ij,m}$ (resp. $\underbar Y^{ij,n}$) are defined in (\ref{penalizedscheme}) (resp. (\ref{second-penalizedscheme})). 

So let $p\geq 1$ and $(Y^{ij,m,p},Z^{ij,m,p},K^{ij,m,p })$  be the processes defined recursively by: $\forall \ig$ and $m,p\ge 0$, \be\label{penalizedschemebis}
\left\{\begin{array}{l}
Y^{ij,m,p+1}\in \cS^{2}_t, \,\,Z^{ij,m,p+1} \in \cH^{2,d}_t \mbox{
and
} K^{ij,m,p+1 } \in \cA^{2}_{t,i}\,\,;\\\\
Y^{ij,m,p+1}_s=h^{ij}(X^{t,x}_T)+\int_{s}^{T}{f}^{ij}(r,X^{t,x}_r,[
(Y^{kl,m,p}_r)_{(k,l)\in \Gamma^{-(i,j)}},Y^{ij,m,p+1}_r],Z^{ij,m,p+1}_r)dr+\int_s^Td K^{ij,m,p+1}_r
\\\qq -m
\int_s^T(Y^{ij,m,p+1}_r-\min_{l\in (\Gamma^{2})^{-j}}\{
Y^{il,m,p}_r+\bar {g}_{jl}(r,X^{t,x}_r)\})^+dr
-\int_s^T Z^{ij,m,p+1}_rdB_r,\,t\le s\leq T;\\\\
\displaystyle{Y^{ij,m,p+1}_{s}\geq \max_{k\in
(\Gamma^{1})^{-i}}\{
Y^{kj,m,p}_s-\underline{g}_{ik}(s,X^{t,x}_s)\}},\,t\le s\leq T;
\\\\
\int_t^T(Y^{ij,m,p+1}_s-\max_{k\in (\Gamma^{1})^{-i}}\{
Y^{kj,m,p}_s-\underline{g}_{ik}(s,X^{t,x}_s)\})d
K^{ij,m,p+1}_s=0\end{array}\right.\ee
where for any $\ig$, we initialize the scheme as follows: $Y^{ij,m,0}=\bar Y^{ij}$.  First note that for any fixed $\ig$ and $p\ge 0$, 
$(Y^{ij,m,p+1},Z^{ij,m,p+1}, K^{ij,m,p+1})$ exists since it solves a standard reflected BSDE. Next and by comparison of solutions of reflected BSDEs \cite{karoui}, we have for any $\ig$, $Y^{ij,m,1}\ge \bar Y^{ij}$. Indeed, this holds true since $(\bar Y^{ij},\bar Z^{ij},\bar K^{ij,\pm})_{\ig}$ verifies (\ref{dreflect-bsde}) and $m
(\bar Y^{ij}_s-\min_{l\in (\Gamma^{2})^{-j}}\{
\bar Y^{il}_s+\bar {g}_{jl}(s,X^{t,x}_s)\})^+=0$ for any $s\in [t,T]$. Thus, to conclude, it is enough to apply It\^o-Meyer's formula
with $((\bar Y^{ij}-Y^{ij,m,1})^+)^2$ to get $((\bar Y^{ij}_s-Y_s^{ij,m,1})^+)^2\equiv 0$, $s \in [t,T]$, which implies the desired result.\\
\ms  Next,
since for any $(q,r)\in \Gamma^{-(i,j)}$, the mapping \\
$ y^{qr}\in \R\mapsto f^{ij}(t,x,[(y^{kl})_{((k,l)\in \Gamma^{-(i,j)}},y],z)-m(y-\min_{q\in \gamj}(y^{iq}+\bar g_{jq}(t,x)))^+$ is non-decreasing when
the other variables are frozen, then one can show by induction that for any $\ig$,  $Y^{ij,m,p}\leq Y^{ij,m,p+1}$ (for $p=0$
this inequality holds true). Referring then to equality (19), Section 4 in \cite{Djehichehampopier}, the family $ (Y^{ij,m,p})_{i \in \Gamma^1}$
 is identified with the value of a standard switching problem with triple of data $(f^{ij, m}, h^{ij} ,(\udl{g}_{ik})_{k \neq i})$ when only $p$ switchings at most are permitted, i.e., 
 \be \label{gamerep-yijmp} \ba{ll}
  Y_t^{ij,m, p} =& \; \displaystyle{ \esssup_{ \{ \delta=(\sigma_k,\;\xi_k)_{k\ge 0} \in {\cal A}^i_{\pi_1}(t), \sigma_{p+1}=T\} }
  \mathbb{E}[ \int_{t}^{\sigma_p \wedge T}\sum_{k=1,p} f^{\xi_{k-1}j, m}(s,X_s)\mathbf{1}_{ \sigma_{k-1} \le s <\sigma_{k}}ds} \\
     & \quad \quad \quad\quad 
     \displaystyle{ - \sum_{k=1,p}\udl{g}_{\xi_{k-1},\xi_{k}}(\sigma_k,X_{\sigma_k})\mathbf{1}_{\{\sigma_k <T\}}}   
\displaystyle{ +\bar  Y_{\sigma_p}^{\xi_{p-1} j}\mathbf{1}_{\{\sigma_p < T\}} +h^{a_T j}(X_T)\mathbf{1}_{\{\sigma_p  =T\}}| \mathcal{F}_t].}\\
  \ea  \ee
 Taking now into account both the representation (\ref{gamerep-yijmp}), Assumption (H5)-(ii) satisfied by $f^{ij}$ and the fact that the penalty costs are non negative, 
 one obtains the existence of a stochastic process  $W^{ij}$ such that for any $(m,p)$ it holds
$$ Y^{ij,m,p}   \le \bar Y^{ij, m}\le W^{ij}.$$ 
 Actually, it is enough to take
$$\ba{l}
W^{ij}_s=\E[\int_s^TC(1+|X^{t,x}_r|^p)dr+\sum_{(k,l)\in \Gamma}\sup_{s\in [t,T]}|\bar Y^{kl}_s||\cF_s],\,\,s\in [t,T]\ea
$$where $C$ is a constant appropriately choosen. 
Next in the same way as in (\cite{Hutang}, Theorem 2.1), one can show that the sequence of processes $((Y^{ij,m,p},Z^{ij,m,p},K^{ij,m,p})_{\ij})_{p\ge 0}$ converges to 
$(Y^{ij,m},Z^{ij,m},K^{ij,m})_{\ij}$ which implies that $\bar Y^{ij,m}\geq \bar Y^{ij}$ for any $\ij$. But, in a symmetric way, one can show that $\bar Y^{ij}\ge \underbar Y^{ij,n}$ for any $\ig$. Take now the limit w.r.t to $m$ and $n$ in the previous inequalities to obtain $Y^{ij}\ge \bar Y^{ij}\ge Y^{ij}$ since $Y^{ij}=\lim_m \bar Y^{ij,m}=\lim_n \underbar Y^{ij,n}$. Thus for any $\ig$, $Y^{ij}=\bar Y^{ij}$. Next, comparing the martingale parts in the equation (\ref{dreflect-bsdeeqxx}) solved by $Y^{ij}$ and $\bar Y^{ij}$, we deduce that $Z^{ij}=\bar Z^{ij}$ and finally $K^{ij,+}-K^{ij,-}= \bar K^{ij,+}-\bar K^{ij,-}$ for any $\ig$. Thus the solution of (\ref{dreflect-bsdeeqxx}) is unique. \\

\begin{remark} Since the interconnected barriers $U^{ij}_s(\vec{Y}):=\min_{l\neq j} \{Y^{il}_s+\bar g_{jl}(s,X_s^{t,x})\}$ and \\
$L^{ij}_s(\vec{Y}):=\max_{k\neq i} \{Y^{kj}_s+\underline g_{kj}(s,X_s^{t,x})\}$ are not completely separated, then one cannot infer that the processes  
$K^{ij,+}$ and $K^{ij,-}$ of (\ref{dreflect-bsdeeqxx}) are unique. But if we moreover require that $dK^{ij,+}$ and $dK^{ij,-}$ are singular then they are actually unique. Indeed, $dK^{ij,+}-dK^{ij,-}$ is a signed measure which has a unique decomposition 
into $d\lambda^{ij,+}- d\lambda^{ij,-}$, i.e.,
$dK^{ij,+}-dK^{ij,-}=d\lambda^{ij,+}- d\lambda^{ij,-}$ where 
$d\lambda^{ij,+}$ and $ d\lambda^{ij,-}$ are non-negative singular measures. Therefore 
$dK^{ij,+}+d\lambda^{ij,-}=d\lambda^{ij,+}+dK^{ij,-}$. Then $d\lambda^{ij,+}<\!\!<dK^{ij,+}$ and 
$d\lambda^{ij,-}<\!\!<dK^{ij,-}$ which implies that 
$d\lambda^{ij,+}_s=a_s^{ij,+}dK^{ij,+}_s$
and  
$d\lambda^{ij,-}_s=a^{ij,-}_sdK^{ij,-}_s$. It follows that $(Y^{ij,+}_s-L^{ij}_s(\vec{Y}))
d\lambda^{ij,+}_s=(Y^{ij}_s-L^{ij}_s(\vec{Y}))a_s^{ij,+}dK^{ij,+}_s=0$ and similarly 
$(Y^{ij,-}_s-U^{ij}_s(\vec{Y}))
d\lambda^{ij,-}_s=(Y^{ij}_s-U^{ij}_s(\vec{Y}))a_s^{ij,-}dK^{ij,-}_s=0$.
\end{remark}
\noindent {\bf (A2): PDEs with bilateral obstacles}
\medskip

Let $(t,x)\in \esp$ and $(X^{t,x}_s)_{s\leq T}$ be the solution of the standard SDE given in (\ref{sdex}) where the functions $b$ and $\sigma$ satisfy Assumption (H0). Let us now consider the following functions:
$$
\begin{array}{ll}
g:& x\in \R^k\longmapsto g(x)\,\R\\
f:& (t,x,y,z)\in [0,T]\times
\R^{k+1+d}\longmapsto f(t,x,y,z)\in \R
\\
H:& (t,x)\in [0,T]\times \R^k\longmapsto H(t,x)\in \R\\
L:&(t,x)\in 
[0,T]\times \R^k\longmapsto L(t,x)\in \R\ea$$ We assume that all those functions are continuous and satisfy the following assumptions (A1)-(A2).
\ms

\no {\bf (A1)}: $\forall \,\,t\in [0,T]$, $x\in \R^k$, $y,y'\in \R$, $z,z'\in
\R^d$, 
$$\left\{
    \begin{array}{ll}
\mbox{(i) } |g(x)|+|f(t,x,0,0)|+|H(t,x)|+|L(t,x)|\leq C(1+|x|^p),\\
\mbox{(ii) } |f(t,x,y,z)-f(t,x,y',z')|\leq C(|y-y'|+|z-z'|),\\
\mbox{(iii) }L(t,x)\leq H(t,x) \mbox{ and }L(T,x)\leq g(x)\leq H(T,x),
    \end{array}
    \right.$$
where $C$ and $p$ are some positive constants.\ms

\noindent {\bf (A2)}: For each $R>0$,
there is a function $\Phi_R$ from $\R^+$ to $\R^+$ satisfying $\Phi_R(s)\rw 0$ as $s\rw 0$ and such that \be \label{eqcontuniff}|f(t,x,y,z)-f(t,x',y,z)|\leq
\Phi_R((1+|z|)|x-x'|)\ee for all $t\in (0,T)$, $|x|$, $|x'|$,
$|y|\leq R$ and $z\in \R^d$.
\ms

The follownig lemma gives a link between Assumptions (A2) and (H1)-(i).

\begin{lem}If the function $(t,x)\in \espo \mapsto f(t,x,y,z)\in \R$ is continuous, uniformly w.r.t. $(y,z)$, i.e., for any $(t,x)\in \espo$, for any 
$\eps >0$ there exists $\eta_{t,x,\eps}>0$ such that if $|(t',x')-(t,x)|<\eta_{t,x,\eps}$ then 
$|f(t',x',y,z)-f(t,x,y,z)|<\eps$, then assumption (A2) is satisfied with $\Phi_R$, for any $R\geq 0$, given by
$$
\forall \g \geq 0, \Phi_R(\g):=\sup_{|t-t'|+|x-x'|\leq \g, t,t'\in [0,T], x,x'\in B'(0,R),|y|\leq R, z\in \R^d}|f(t,x,y,z)-f(t',x',y,z)|
$$
where $B'(0,R):=\{x\in \R^k,|x|\leq R\}$.
\end{lem}
\no $\pr$: Let $R>0$ and $(t,x)\in [0,T]\times B'(0,R)$. Let $\eps>0$. By definition, there exists 
$\eta_{t,x,\eps}>0$ such that if $(t',x')\in B((t,x),\eta_{t,x,\eps})$ (the open ball in $\espo$ with center $(t,x)$ and radius $\eta_{t,x,\eps}$) then 
$$|f(t',x',y,z)-f(t,x,y,z)|<\eps.$$
As $$[0,T]\times B'(0,R)\subset \bigcup_{(t,x)\in [0,T]\times B'(0,R)}B((t,x),\frac{\eta_{t,x,\eps}}{2})
$$
then by compacity one can find finitely many points $(t_1,x_1), ...., (t_m,x_m)$ such that 
$$[0,T]\times B'(0,R)\subset \bigcup_{i=1,m}B((t_i,x_i),\frac{\eta_{t_i,x_i,\eps}}{2}).$$
\ms

\no (i) First note that the function $\gamma \in \R^+\mapsto \Phi_R(\g)$ is non-decreasing and $\Phi_R(0)=0$. Next let us set $\eta=\min_{i=1,m}\frac{\eta_{t_i,x_i,\eps}}{3}$. Then $\eta>0$ and we claim that $\Phi_R(\eta)\in \R^+$. 

Indeed let $z\in \R^d$, $|y|\leq R$ and $(t,x)$, $(t',x')$ elements of $[0,T]\times B'(0,R)$ such that 
$|t-t'|+|x-x'|\leq \eta$. Then there exists $i\in\{1,...,m\}$ such that $(t,x)\in B((t_i,x_i),\frac{\eta_{t_i,x_i,\eps}}{2})$. It follows that 
$$
|(t',x')-(t_i,x_i)|\leq |(t',x')-(t,x)|+|(t,x)-(t_i,x_i)|\leq \eta +\frac{\eta_{t_i,x_i,\eps}}{2}<\eta_{t_i,x_i,\eps}
$$
which implies that $(t',x')$ belongs also to $B((t_i,x_i),\eta_{t_i,x_i,\eps})$. Then by continuity we have  
$$
|f(t,x,y,z)-f(t',x',y,z)|\leq |f(t,x,y,z)-f(t_i,x_i,y,z)|+|f(t_i,x_i,y,z)-f(t',x',y,z)|\leq 2\eps. 
$$
Taking the supremum to obtain that $\Phi_R(\eta)\leq 2\eps$. 
\ms

\no (ii) As $\Phi_R$ is non-decreasing then $\Phi_R(\gamma)\leq 2\eps $ for any $\gamma\leq \eta$. Note that this property implies also that $\Phi_R(\gamma)\rw 0=\Phi_R(0)$ as $\gamma \rw 0$. 
\ms

\no (iii) For any $\g_1$ and $\g_2$ in $\R^+$, $\Phi_R(\g_1+\g_2)\leq  \Phi_R(\g_1)+\Phi_R(\g_2)$.

Indeed let $z\in \R^d$, $|y|\leq R$ and $(t,x)$, $(t',x')$ elements of $[0,T]\times B'(0,R)$ such that 
$|t-t'|+|x-x'|\leq \g_1+\g_2$. Then there exists $(\bar t,\bar x)\in [0,T]\times B'(0,R)$ such that 
$
|(t,x)-(\bar t,\bar x)|\leq \g_1$ and $|(\bar t,\bar x)-(t',x')|\leq \g_2$. Therefore 
$$
|f(t,x,y,z)-f(t',x',y,z)|\leq |f(t,x,y,z)-f(\bar t,\bar x,y,z)|+|f(\bar t,\bar x,y,z)-f(t',x',y,z)|\leq  \Phi_R(\g_1)+\Phi_R(\g_2)
$$
which implies that $\Phi_R(\g_1+\g_2)\leq  \Phi_R(\g_1)+\Phi_R(\g_2)$.
\ms 

\no (iv) For any $\gamma \in \R^+$, $\Phi_R(\gamma)\in \R^+$ and $\Phi_R(.)$ verifies (\ref{eqcontuniff}). 

Indeed by induction and (iii), for any $\gamma \in \R^+$ and $n\ge 1$, $\Phi_R(n\g)\leq n\Phi_R(\gamma)$. On the other hand, by (ii) and (iii), for any $\g\in \R^+$ one can find an integer $n$ such that $\Phi_R(\g)\leq n\Phi_R(\eta)$ which implies that $\Phi_R(\gamma)\in \R$. Finally 
(\ref{eqcontuniff}) is obviously satisfied since  $\Phi_R(.)$ is non-decreasing and $$
|f(t,x,y,z)-f(t,x',y,z)|\leq
\Phi_R(|x-x'|)$$\ for all $t\in (0,T)$, $|x|$, $|x'|$,
$|y|\leq R$ and $z\in \R^d$.  
\begin{rem}
This property of continuity of the function $(t,x)\in \espo \mapsto f(t,x,y,z)\in \R$, uniformly w.r.t. $(y,z)$, is needed to get uniqueness of the viscosity solution of the PDE (\ref{appendixu}). However one can obtain uniqueness of 
(\ref{appendixu}) with a substantially weaker condition than the previous one. \qed
\end{rem} 
Now for $(t,x)\in \esp$, let $Y^{t,x}:=(Y^{t,x}_s)_{s\in [0,T]}$ be the local solution of the BSDE associated with the quadruple $(f(s,X^{t,x}_s,y,z), g(X^{t,x}_T), L(s,X^{t,x}_s), H(s,X^{t,x}_s))$ (see \cite{Hamhassani} for more details) and let us set $u(t,x)=Y^{t,x}_t$. Then $u(t,x)$ is an $\R$-valued deterministic function of $(t,x)$ which is morever of polynomial growth and continuous. Furthermore it satisfies: 
\be\lb{lienyu}\forall \,\,\tx,\,\,\forall s\in [t,T],\,\,
Y^{t,x}_s=u(s,X^{t,x}_s).
\ee

Let us now consider the following PDE with two obstacles of min-max type whose solutions will be considered in viscosity sense:  \be\label{appendixu}\left\{
    \begin{array}{ll}
\min\Big \{v(t,x)-L(t,x)\,;\,\max\Big[v(t,x)-H(t,x);\\\qq\qq\qq
 -\partial_t v(t,x)-\cL^X v(t,x)-f(t,x,v(t,x),\sigma(t,x)^\top D_x v(t,x))\Big]\Big \}=0\,;\\
     v(T,x)=g(x).
    \end{array}
    \right.\ee
For the definition of the viscosity solution of equation (\ref{appendixu}), which is standard, we refer the reader to \cite{Hamhassani}. 
\ms

The link between PDE (\ref{appendixu}) and the process $Y^{t,x}$ through the function $u$ defined in (\ref{lienyu}) is: 
\begin{thm} \cite{Hamhassani}: \label{existence} Under (H0),(A1) and (A2) we have: 
\ms

(i) The function $u$ is the unique continuous viscosity solution of (\ref{appendixu})
with polynomial growth ;
\ms

(ii) The function $u$ is also a unique continuous viscosity solution, in the class $\pl$, of the following max-min problem:
\be\label{appendixu3}\left\{
    \begin{array}{ll}
\max\Big \{v(t,x)-H(t,x)\,\,;\,\,\min\Big[v(t,x)-L(t,x);\\\qq\qq\qq
 -\partial_t v(t,x)-\cL^X v(t,x)-f(t,x,v(t,x),\sigma(t,x)^\top D_x v(t,x))\Big]\Big \}=0;\\
     v(T,x)=g(x).
    \end{array}
    \right.\ee
\end{thm}
The proof of (i) is similar to the one given in \cite{Hamhassani}. However, we should point out that in \cite{Hamhassani}, the barriers $L$ and $H$ are assumed to be completely separated (i.e. $L<H$) while in our framework they only satisfy $L\leq U$. This fact is irrelevant and does not rise a major issue. As for (ii), the construction of the function $u$ (see \cite{Hamhassani}) implies that $w=-u$ is the unique viscosity solution in the class $\pl$ of the following system:
\be\label{appendix}\left\{
    \begin{array}{ll}
\min\Big \{w(t,x)+H(t,x),\max\Big[w(t,x)+L(t,x),\\\qq\qq\qq
 -\partial_t w(t,x)-\cL w(t,x)+f(t,x,-w(t,x),-\sigma(t,x)^\top D_x w(t,x))\Big]\Big \}=0;\\
     w(T,x)=-g(x).
    \end{array}
    \right.\ee
Thus $-w=u$ is the unique solution in the class $\pl$ of system (\ref{appendixu3}) (see e.g. \cite{barles}, pp.18). \qed

\noindent {\bf References}:

\end{document}